\documentclass[letterpaper]{article}
\usepackage[margin=1in]{geometry}
\usepackage{float}
\usepackage{multirow} 
\usepackage{bm}
\usepackage{amsmath,amsthm,amsfonts,amssymb}
\usepackage{enumerate,xcolor}
\usepackage{graphicx}
\usepackage{subcaption}
\usepackage[authoryear,round]{natbib}
\usepackage[ruled,vlined]{algorithm2e}
\usepackage{algorithmic}
\usepackage{booktabs,tabularx}

\PassOptionsToPackage{hyphens}{url}
\usepackage{url} 
\usepackage{hyperref}
\hypersetup{
  pdftitle  = {Scaling limit of heavy-tailed nearly unstable cumulative INAR(\texorpdfstring{$\infty$}{infty}) processes and rough fractional diffusions},
  pdfauthor = {Yingli Wang and Chunhao Cai and Ping He and Qinghua Wang},
  colorlinks = true,
  linkcolor  = [rgb]{0.10,0.20,0.55},
  citecolor  = [rgb]{0.10,0.20,0.55},
  urlcolor   = [rgb]{0.10,0.20,0.55}
}

\usepackage{booktabs}
\usepackage{array}

\allowdisplaybreaks

\numberwithin{equation}{section}
\newtheorem{theorem}{Theorem}[section]
\newtheorem{lemma}[theorem]{Lemma}
\newtheorem{proposition}[theorem]{Proposition}
\newtheorem{corollary}[theorem]{Corollary}

\newtheorem{assumption}[theorem]{Assumption}

\theoremstyle{definition}
\newtheorem{definition}[theorem]{Definition}

\theoremstyle{remark}
\newtheorem{remark}[theorem]{Remark}

\newcommand{\norm}[1]{\left\lVert#1\right\rVert}

\begin{document}

\title{Scaling limit of heavy-tailed nearly unstable cumulative INAR($\infty$) processes and rough fractional diffusions}

\author{%
Yingli Wang\thanks{Corresponding author. School of Mathematics, Shanghai University of Finance and Economics, Guoding Road, Shanghai 200433, Shanghai, China; 2022310119@163.sufe.edu.cn.}%
\and
Chunhao Cai\thanks{School of Mathematics (Zhuhai), Sun Yat-Sen University, Tangjiawan, Zhuhai Campus, Zhuhai 519082, Guangdong, China; caichh9@mail.sysu.edu.cn.}%
\and
Ping He\thanks{School of Mathematics, Shanghai University of Finance and Economics, Guoding Road, Shanghai 200433, Shanghai, China; pinghe@mail.shufe.edu.cn.}%
\and
Qinghua Wang\thanks{School of Mathematics, Shanghai University of Finance and Economics, Guoding Road, Shanghai 200433, Shanghai, China; wangqinghua@mail.shufe.edu.cn.}%
}

\date{\today}

\maketitle

\begin{abstract}
In this paper, we investigate the scaling limit of heavy-tailed nearly unstable cumulative INAR($\infty$) processes. These processes exhibit a power-law tail of the form $n^{-(1+\alpha)}$ for $\alpha \in (\frac{1}{2}, 1)$, and the $\ell^1$ norm of the kernel vector converges to 1. We demonstrate that the discrete-time scaling limit retains a long-memory property and can be viewed as an integrated fractional Cox-Ingersoll-Ross process. Moreover, we present an efficient method for simulating the fractional Cox-Ingersoll-Ross process. The simulation and Goodness-of-Fit Test code are available at \url{https://github.com/gagawjbytw/INAR-rough-Heston}.
\end{abstract}

\noindent%

\noindent%


\section{Introduction}
Let $(\epsilon_n)_{n\in\mathbb{Z}}$ be i.i.d. \text{Pois}($\nu$) random variables. Let $(\eta_k)_{k\ge1}$ be a sequence of positive numbers in $\ell^1$. For $n \in \mathbb{Z}$, $k \ge 1$, and $l \ge 1$, let $\xi_l^{(n,k)}$ be mutually independent random variables with a \text{Pois}($\eta_k$) distribution, which are also independent of the sequence $(\epsilon_n)_{n\in\mathbb{Z}}$.

An Integer-Valued Autoregressive process of infinite order, abbreviated as INAR($\infty$), refers to a type of integer-valued time series model denoted by $(X_n)_{n \in \mathbb{Z}}$. This model constitutes an infinite-order extension of the classic INAR($p$) processes. The process is formally defined by the following system of stochastic difference equations:
\[
  \epsilon_n = X_n - \sum_{k=1}^\infty \sum_{l=1}^{X_{n-k}} \xi_l^{(n,k)}, \quad n \in \mathbb{Z}.
\]
To establish a more compact notation, we introduce the \textbf{reproduction operator}, denoted by $\star$. For a reproduction parameter $\iota \ge 0$, we define the \textbf{offspring sequence} $(\xi_n^{(\iota)})_{n \ge 1}$ as a sequence of i.i.d. random variables, where each element $\xi_n^{(\iota)}$ is an \textbf{offspring variable} following a Poisson distribution with mean $\iota$.

The reproduction operator is then defined for any non-negative integer-valued random variable $Y$ as the sum of its offspring:
\[
  {\iota \star Y} := \sum_{n=1}^Y \xi_n^{(\iota)},
\]
provided that $Y$ is independent of the offspring sequence $(\xi_n^{(\iota)})$. {Note that when $\xi_n^{(\iota)}\sim \mathrm{Bernoulli}(\alpha)$, the operator $\iota \star Y$ coincides with the classical binomial thinning operator $\iota \circ Y$ used in the INAR($p$) literature.}
Using this operator, the system of stochastic difference equations can be compactly rewritten as:
\[
  \epsilon_n = X_n - \sum_{k=1}^\infty {\eta_k \star X_{n-k}}.
\]
Finally, we {call}
$\nu$ the \textbf{immigration parameter}, $(\epsilon_n)$ the \textbf{immigration sequence}, and $\eta_k$ the \textbf{reproduction coefficient} for each non-negative integer $k$.

A cumulative INAR($\infty$) process (also known as a discrete Hawkes process), which starts from time $1$ and is denoted by $(N_n)_{n\ge1}$, is defined by
\begin{equation}
  N_n = \sum_{s=1}^n X_s.
\end{equation}
Hawkes processes, introduced by \cite{hawkes1971spectra}, are continuous-time self-exciting point processes widely used in various fields. Discrete-time analogs, such as cumulative INAR($\infty$) processes, provide similar modeling capabilities, particularly for count data observed at fixed time intervals. However, the scaling limits of heavy-tailed, {nearly unstable (i.e., near-critical)} cumulative INAR($\infty$) processes differ considerably from those of continuous-time counterparts, which warrants a dedicated investigation. {Here ``nearly unstable'' means that the total reproduction mean is close to the critical value $1$ (see \eqref{eq:nearly_unstable}).}
Under certain conditions, the Poisson autoregressive process can be viewed as an INAR($\infty$) process with Poisson offspring.
For a comprehensive discussion on Poisson autoregressive models and their connections to INAR and Hawkes processes, refer to \cite{fokianos2021multivariate}, \cite{huang2023nonlinear} and \cite{kirchner2016hawkes}. 
From Proposition 3 in \cite{kirchner2016hawkes}, we see that if an INAR($\infty$) process $(X_n)_{n \geq 1}$ starts from time $1$, {it can be specified recursively through its conditional mean (or intensity) $(\lambda_n)_{n\ge 1}$ as follows:}
  \begin{equation}\label{process_definition}
  \lambda_n := \nu + \sum_{s=1}^{n-1} \eta_{n-s} X_s, \qquad n\ge 1,
  \end{equation}
  where $\nu > 0$ is the immigration rate, and $(\eta_n)_{n \geq 1} \in \ell^1$ represents the offspring distribution with $\eta_n \geq 0$ for all $n \in \mathbb{N}$. {Let $\mathcal{F}_{n}:=\sigma(X_1,\dots,X_n)$ (with $\mathcal{F}_0$ the trivial $\sigma$-field); when we later consider the $T$-dependent sequence $(X_n^T)$, we use the same notation $\mathcal{F}_n:=\sigma(X_1^T,\dots,X_n^T)$ for the filtration it generates.} Then $\lambda_1=\nu$ and, conditionally on $\mathcal{F}_{n-1}$, the count $X_n$ is Poisson with parameter $\lambda_n$, i.e.,
  \[
  X_n \mid \mathcal{F}_{n-1} \sim \text{Pois}(\lambda_n).
  \]

Scaling limits for self-exciting processes, a topic of significant theoretical and applied interest, have been extensively studied. For example, consider a Hawkes process $(N_t)_{t\ge0}$ with intensity function $h$. A common goal is to characterize the asymptotic behavior of the rescaled process 
\[
  (\beta_T N_{tT})_{t\in[0,1]},
\]
as $T\rightarrow\infty$, where $\beta_T$ is an appropriate normalizing factor. In \cite{bacry2013some}, it is shown that under the condition
\[
  \norm h_{L^1}\equiv\int_0^\infty h(s)ds<1,
\]
the asymptotic behavior of a Hawkes process is quite similar to that of a Poisson process. Indeed, as $T$ tends to infinity,
\[
  \sup_{t\in[0,1]}\left| \frac{N_{tT}}{T}-\mathbb E\left[ \frac{N_{tT}}{T} \right] \right|\rightarrow 0\ \text{in probability},
\]
and
\[
  \left( \sqrt T\left( \frac{N_{tT}}{T}-\mathbb E\left[ \frac{N_{tT}}{T} \right] \right) \right)_{t\in[0,1]}\rightarrow \sigma(W_t)_{t\in[0,1]}\ \text{in law}
\]
under the Skorokhod topology, with $\sigma$ an explicit constant and $(W_t)$ a Brownian motion. This result has been extended in \cite{zhu2013nonlinear} to the case of nonlinear Hawkes processes. Regarding unstable Hawkes processes, the light-tailed case and the heavy-tailed case are respectively studied in \cite{jaisson2015limit} and \cite{jaisson2016rough}. The light-tailed case corresponds to a CIR limit obtained under the crucial assumption
\[
  \int_0^\infty sh(s)ds<+\infty
\]
and the heavy-tailed case corresponds to a fractional CIR limit obtained under
\[
  \int_0^\infty sh(s)ds=+\infty\ \text{and}\ h(x)\sim \frac{K}{x^{1+\alpha}}\ \text{as}\ x\rightarrow\infty.
\]

In the fields of mathematical finance and econometrics, for example, the limit behavior of self-exciting dynamics is of particular interest. Studies such as \cite{horst2019scaling, horst2022microstructure} examine the microstructure of financial markets and the impact of external shocks on market dynamics. They mathematically demonstrate that multiple small external shocks can trigger endogenous jump cascades in asset returns and stock price volatility.

When it comes to the discrete case, the asymptotic behavior of unstable Integer-Valued Autoregressive processes of order $p$, or INAR($p$), was investigated by \cite{barczy2011asymptotic}. An INAR($p$) process $(X_n)$ is defined by the stochastic difference equation 
\[
    X_n = \sum_{i=1}^p \alpha_i \circ X_{n-i} + \epsilon_n,
\]
where the $\alpha_i \in [0, 1)$ are autoregressive parameters, $\circ$ denotes the binomial thinning operator, and $(\epsilon_n)$ is a sequence of i.i.d. non-negative integer-valued random variables representing an immigration term. The authors proved that under the critical condition where the sum of these autoregressive parameters equals one:
\[
  \sum_{i=1}^p\alpha_i=1,
\]
a sequence of appropriately scaled random step functions formed from the process converges weakly towards a squared Bessel process. {
In the light-tailed (summable-kernel) setting, closely related results are available and suggest a squared-Bessel-type scaling limit for INAR($\infty$).
Indeed, unstable limits are proved for finite-order discrete-time models (INAR/INAP($p$)) in \cite{barczy2011asymptotic}, and the nearly-unstable light-tailed scaling limit is established for continuous-time Hawkes processes in \cite{jaisson2015limit}.
To the best of our knowledge, a dedicated proof for the INAR($\infty$) light-tailed case is not explicitly available in the literature.
Consequently, scaling limits for nearly-unstable INAR($\infty$) models with heavy-tailed reproduction kernels appear to be missing from the literature, and this is the gap we aim to fill.}

In this paper, we aim to provide a discrete-time counterpart to the continuous-time results found in \cite{jaisson2016rough}, focusing on heavy-tailed nearly unstable INAR($\infty$) processes. Specifically, we consider a baseline reproduction kernel $(\eta_n)_{n\ge1}$ with heavy tail
\[
  \eta_n \sim \frac{K}{n^{1+\alpha}} \quad \text{as} \quad n \to \infty,
\]
for some $\alpha \in \left(\frac{1}{2}, 1\right)$ and $K>0$, and we study a sequence of INAR($\infty$) processes $(X_n^T)_{n\in\mathbb{N}}$ as $T\to\infty$.
{Throughout, ``nearly unstable'' means that the total reproduction mean approaches the critical value $1$, namely
\begin{equation}\label{eq:nearly_unstable}
  \|\eta^T\|_1=\sum_{k\ge1}\eta_k^T \uparrow 1 \quad \text{as } T\to\infty.
\end{equation}
In our parametrization, we take $\eta^T=a_T\eta$ with $a_T<1$ and $a_T\uparrow 1$, so that $\|\eta^T\|_1=a_T$.}
{For each $T$, $(X_n^T)_{n\ge1}$ is an INAR($\infty$) process with immigration parameter $\mu^T$ and reproduction coefficients $(\eta_k^T)_{k\ge1}$. We denote by $\lambda_n^T:=\mathbb E[X_n^T\mid\mathcal{F}_{n-1}]$ its conditional mean (intensity), where here $\mathcal{F}_{n-1}=\sigma(X_1^T,\dots,X_{n-1}^T)$.}
The sequence of parameters $(\lambda_n^T)$ is defined by
\begin{equation}\label{eq:lambdaTexpression}
  \lambda_n^T = \mu^T + \sum_{s=1}^{n-1} \eta_{n-s}^T X_s^T, \quad n \in \mathbb{N},
\end{equation}
where $\mu^T$ is a sequence of positive real numbers. As in \cite{jaisson2016rough}, we establish specific assumptions on $\eta^T$ to facilitate our analysis. We prove tightness of the renormalized cumulative INAR$(\infty)$ process $Y^T$ and the associated martingale $Z^T$ defined in \eqref{eq:expressionofy} and \eqref{eq:expressionofz}. Moreover, we identify any subsequential limit $(Z,Y)$ of $(Z^T,Y^T)$ as a continuous martingale $Z$ whose sharp bracket satisfies $[Z,Z]_t = Y_t$, where $Y$ is an integrated fractional Cox--Ingersoll--Ross (CIR) process.

While our work is inspired by and presents a discrete-time counterpart to the celebrated results of \cite{jaisson2016rough}, it is motivated by distinct theoretical and practical considerations that underscore its unique contribution. The transition from a continuous-time Hawkes process to a discrete-time INAR($\infty$) process is far from trivial and offers several key advantages:

\begin{enumerate}
    \item \textbf{Direct Modeling of Intrinsically Discrete Data:} Many real-world phenomena, particularly in finance (e.g., tick-by-tick trade counts, order book events) and social sciences (e.g., daily event counts), are recorded as discrete counts over fixed time intervals. The INAR($\infty$) framework provides a more natural and direct modeling tool for such data, avoiding the potentially lossy abstraction required to fit them into a continuous-time point process framework.

    \item \textbf{Novel Proof Techniques and Theoretical Insights:} The proofs in the discrete setting are not mere translations of their continuous-time counterparts. They require different mathematical machinery, such as discrete renewal theory and specific martingale central limit theorems for discrete processes. This distinct theoretical path reveals unique structural properties of the discrete models, such as the specific long-memory characteristics we highlight.

    \item \textbf{A Computationally Superior Simulation Framework:} Perhaps the most significant practical contribution is that INAR($\infty$) processes are substantially easier to simulate than their continuous-time, singular-kernel Hawkes counterparts. Simulating our model only requires iterative sampling from a Poisson distribution. Our scaling limit theorem thus serves as a crucial bridge: it rigorously validates that this simple discrete simulation converges to the complex continuous process of interest. This directly addresses the question of ``why not just use a Hawkes process?'': The INAR($\infty$) process provides a computationally tractable pathway for its analysis and application.

    \item \textbf{New Applications in Financial Modeling:} This framework opens the door to new applications not explored in the original Hawkes process literature. For instance, as demonstrated in the recent work of \cite{wang2025rough}, a bivariate version of this INAR($\infty$) framework can be shown to converge to the rough Heston model. This provides not only a novel microstructural foundation for rough volatility but also a simulation engine for pricing both European and path-dependent options, a task that remains challenging for other methods.
\end{enumerate}

Our primary result, encapsulated in Theorem \ref{theorem:1}, demonstrates that under heavy-tailed and unstable conditions, the cumulative INAR($\infty$) process converges to a rough fractional diffusion process. This finding is achieved by leveraging discrete renewal theory and martingale central limit theorems.

These results contribute to the theoretical understanding of discrete self-exciting processes by revealing inherent long-memory characteristics under heavy-tailed conditions. Such insights have significant implications for fields such as finance, epidemiology, and network traffic analysis, where capturing long-range dependencies and extreme events is crucial for accurate modeling and forecasting.

\section{Assumptions and Main Results}
\label{sec:assumptionsandmainresults}
\paragraph{Roadmap of Section~\ref{sec:assumptionsandmainresults}.}
We collect here the definitions and statements that will be used throughout the rest of the paper.
Assumption~\ref{assumption:1} specifies the heavy-tailed baseline kernel $(\eta_n)$ and the near-critical scaling $a_T\uparrow 1$,
while Assumption~\ref{assumption:2} fixes the joint scaling of $(a_T,\mu^T)$ that leads to a non-degenerate limit.
We then introduce the rescaled cumulative process $Y^T$ in \eqref{eq:expressionofy}, the rescaled cumulative intensity $\Lambda^T$ in \eqref{eq:LambdaT_def},
and the rescaled martingale $Z^T$ in \eqref{eq:expressionofz}; these are the main objects in the tightness and limit identification arguments.
Lemma~\ref{lemma:laplaceofrho} identifies the limiting Mittag--Leffler kernel via Laplace transforms (used to identify the drift term in the limit),
whereas Lemma~\ref{lemma:tightnessofyandlambda} provides $C$-tightness inputs (used to obtain tightness of $(Z^T,Y^T)$ in Proposition~\ref{proposition:1}).
The remaining lemmas in this section are auxiliary tools: Lemma~\ref{lemma:sumandnotsum} and Lemma~\ref{lemma:discreteKaramataTauberian}
are Tauberian ingredients used in the proof of Lemma~\ref{lemma:laplaceofrho}, and Lemma~\ref{lem:Ak_uniform_bound} provides a uniform pointwise bound
on the resolvent coefficients, used only in the moment/tightness estimates.
Proofs are deferred to Section~\ref{sec:proofoflemmas} and \ref{sec:proofofthemainresults}.

\subsection{Heavy-tail and scaling assumptions}

We normalize the baseline kernel so that $\|\eta\|_1=1$, and we propose two assumptions similar to those in \cite{jaisson2016rough}.

\begin{assumption}\label{assumption:1}
We assume the following conditions hold:
\begin{enumerate}
    \item The baseline offspring distribution $(\eta_n)_{n\ge1}$ is a probability mass function on $\mathbb N$, i.e.
    \[
      \eta_n\ge 0,\qquad \sum_{n=1}^\infty \eta_n=1,
    \]
    and it has a heavy tail in the sense that
    \[
      \eta_n \sim \frac{K}{n^{1+\alpha}}\quad \text{as } n\to\infty,
    \]
    for some $\alpha\in\left(\frac12,1\right)$ and some constant $K>0$.

    \item The sequence of scaling factors $(a_T)_{T\ge0}$ consists of positive numbers such that for all $T$, $a_T < 1$, and it converges to 1 as $T \to \infty$.

    \item For each $T$, the reproduction coefficients are given by $\eta_n^T = a_T \eta_n$ for all $n \ge 1$.
\end{enumerate}
\end{assumption}

{Assumption~\ref{assumption:1}(i) can equivalently be expressed in terms of the tail sum of $(\eta_n)$. To this end, we state the following discrete Karamata-type lemma, which can be derived from the integral versions in \cite[Prop.~1.5.8 and Prop.~1.5.10]{bingham1989regular} via a standard sum--integral comparison.}

\begin{lemma}\label{lemma:sumandnotsum}
  Let $l(t)$ be a {positive} slowly varying function that is locally bounded on $[X_0, \infty)$ for some $X_0 > 0$, and let $\beta\in\mathbb{R}$. For any fixed integer $X \ge X_0$, the following holds as integers $x \to \infty$:
  \begin{enumerate}
    \item {If $\beta > -1$, then}
    \[
      \sum_{t=X}^x t^\beta l(t) \sim \frac{1}{\beta+1}x^{\beta+1}l(x).
    \]
    \item {If $\beta < -1$, then}
    \[
      {\sum_{t=x+1}^\infty t^\beta l(t) \sim -\frac{1}{\beta+1}x^{\beta+1}l(x).}
    \]
  \end{enumerate}
\end{lemma}

{
\begin{proof}
    The proof is provided in Section~\ref{subsec:proofofsumandnotsum}.
\end{proof}
}

{Applying Lemma~\ref{lemma:sumandnotsum} with $\beta=-1-\alpha$ and $l(\cdot)\equiv K$ yields}
\[
  {\sum_{k=n+1}^\infty \eta_k \sim \frac{K}{\alpha}\,n^{-\alpha}\quad\text{as}\quad n\to\infty,}
\]
{and hence}
\[
  \lim_{n\to\infty} \alpha n^\alpha \sum_{k=n+1}^\infty \eta_k = K.
\]

We now specify the asymptotic regime for $(a_T)$ and $(\mu^T)$.
Define the constant
\[
  \delta \;:=\; K\,\frac{\Gamma(1-\alpha)}{\alpha},
\]
which will naturally arise from the small-$z$ expansion of the Laplace transform of the heavy tail (see Lemma~\ref{lemma:laplaceofrho}).

Recall from \eqref{eq:lambdaTexpression} that the (conditional) intensity sequence $(\lambda_n^T)_{n\ge1}$ is given by
\[
  \lambda_n^T \;=\; \mu^T + \sum_{s=1}^{n-1}\eta_{n-s}^T X_s^T,\qquad n\ge1,
\]
where $\mu^T>0$ is the (possibly $T$-dependent) immigration parameter and $(\eta_k^T)_{k\ge1}\in\ell^1$ are reproduction coefficients with $\eta_k^T\ge0$.
{We define the filtration $\mathcal{F}_n:=\sigma(X_1^T,\dots,X_n^T)$ (with $\mathcal{F}_0$ trivial). Then, by construction, $\lambda_n^T=\mathbb{E}[X_n^T\mid\mathcal{F}_{n-1}]$ and
\[
  X_n^T\mid\mathcal{F}_{n-1}\sim\mathrm{Pois}(\lambda_n^T),\qquad n\ge1,
\]
so $(X_n^T-\lambda_n^T)_{n\ge1}$ is a martingale difference sequence.}
We also define the counting (cumulative) process $N_n^T:=\sum_{s=1}^n X_s^T$.

\begin{assumption}\label{assumption:2}
  There exist two positive constants $\lambda,\nu^*$ such that
  \[
    \lim_{T\rightarrow\infty}T^\alpha(1-a_T)=\lambda\delta,
  \]
  and
  \[
    \lim_{T\rightarrow\infty}T^{1-\alpha}\mu^T=\nu^*\delta^{-1},
  \]
  where $\delta=K\frac{\Gamma(1-\alpha)}{\alpha}$. The sequence $(\mu^T)$ is {introduced}
  in \eqref{eq:lambdaTexpression}.
\end{assumption}

\paragraph{Rescaled processes and martingales.}
Define the compensated martingale
\[
  M_n^T:=\sum_{s=1}^n \left(X_s^T-\lambda_s^T\right)=N_n^T-\sum_{s=1}^n\lambda_s^T,\qquad n\ge1.
\]
{Indeed, since $\lambda_s^T=\mathbb{E}[X_s^T\mid\mathcal{F}_{s-1}]$, we have
$\mathbb{E}[X_s^T-\lambda_s^T\mid\mathcal{F}_{s-1}]=0$, hence $(M_n^T)_{n\ge0}$ is an $(\mathcal{F}_n)$-martingale.}

We also introduce the renewal (resolvent) sequence $(A_n^T)_{n\ge1}\in\ell^1$:
\begin{equation}\label{ATexpression}
  A_n^T:=\sum_{k=1}^\infty(\eta^{T*k})_n,
\end{equation}
where $*$ denotes the discrete convolution
\[
  (q*m)(n)=\sum_{s=1}^{n-1}q_sm_{n-s},
\]
and $\eta^{*k+1}=\eta*\eta^{*k}$. This is well-defined since $\|\eta^T\|_1=a_T<1$.

\medskip

We consider, for $t\in[0,1]$, the rescaled cumulative process
\begin{equation}\label{eq:expressionofy}
  Y_t^T=\frac{1-a_T}{T^\alpha \nu^*\delta^{-1}}\sum_{s=1}^{\lfloor Tt \rfloor}X_s^T
  =\frac{1-a_T}{T^\alpha \nu^*\delta^{-1}}N_{\lfloor Tt \rfloor}^T,
\end{equation}
and its rescaled cumulative intensity
\begin{equation}\label{eq:LambdaT_def}
  \Lambda_t^T=\frac{1-a_T}{T^\alpha \nu^*\delta^{-1}}\sum_{s=1}^{\lfloor Tt \rfloor}\lambda_s^T.
\end{equation}
The associated (rescaled) martingale is
\begin{equation}\label{eq:expressionofz}
  Z_t^T:=\sqrt{\frac{T^\alpha \nu^*\delta^{-1}}{1-a_T}}(Y_t^T-\Lambda_t^T).
\end{equation}
{
To justify the terminology, let $\mathcal F_n:=\sigma(X_1^T,\dots,X_n^T)$. Since
$\lambda_n^T=\mathbb E[X_n^T\mid\mathcal F_{n-1}]$, the sequence $(X_n^T-\lambda_n^T)_{n\ge1}$
is a martingale difference sequence, and hence
\[
M_n^T:=\sum_{s=1}^n (X_s^T-\lambda_s^T)
\]
is an $(\mathcal F_n)$-martingale. Therefore,
\[
  Y_t^T-\Lambda_t^T=\frac{1-a_T}{T^\alpha \nu^*\delta^{-1}}\,M_{\lfloor Tt\rfloor}^T
\]
is a martingale in $t$, and so $(Z_t^T)_{t\in[0,1]}$ is a deterministically rescaled martingale.}

\paragraph{Mittag--Leffler kernel.}
Recall the function $f^{\alpha,\lambda}$ from Section 3.1 in \cite{jaisson2016rough}, defined for {$x>0$} by
\begin{equation}\label{eq:expressionoff}
  f^{\alpha,\lambda}(x)=\lambda x^{\alpha-1}E_{\alpha,\alpha}(-\lambda x^\alpha),
\end{equation}
and for {$t\ge0$} by
\begin{equation}\label{eq:Falphalambda}
  F^{\alpha,\lambda}(t)=\int_0^t f^{\alpha,\lambda}(x)\,dx.
\end{equation}
{It is known that $f^{\alpha,\lambda}$ is a probability density on $(0,\infty)$ and $F^{\alpha,\lambda}$ is the corresponding distribution function (the Mittag--Leffler distribution).}
Here $E_{\alpha,\beta}$ is the Mittag--Leffler function
\[
  E_{\alpha,\beta}(z)=\sum_{n=0}^\infty \frac{z^n}{\Gamma(\alpha n+\beta)}.
\]
{We will repeatedly use its Laplace transform: for $z\ge 0$,}
{\begin{equation}\label{eq:laplace_mittagleffler}
  \int_0^\infty e^{-zx} f^{\alpha,\lambda}(x)\,dx = \frac{\lambda}{\lambda+z^\alpha},
\end{equation}}
{equivalently, $\int_0^\infty e^{-zx}\,dF^{\alpha,\lambda}(x)=\lambda/(\lambda+z^\alpha)$.}

\subsection{A key discrete Laplace limit and tightness}

We first state the discrete analog of Lemma 4.3 in \cite{jaisson2016rough}.

\begin{lemma}\label{lemma:laplaceofrho}
  Define a sequence of functions $(\rho^T)_{T\ge0}$ on $[0, \infty)$ by
  \begin{equation}\label{eq:rho_T_def}
    \rho^T(x) := T\frac{A_{\lfloor xT \rfloor}^T}{\norm {A^T}_1}.
  \end{equation}
  {Let $P^T$ be the probability measure on $[0,\infty)$ with density $\rho^T$, and let $F^T(t):=\int_0^t\rho^T(x)\,dx$ be its distribution function.}
  {Then, as $T\to\infty$, $P^T$ converges weakly to the probability measure $P$ on $[0,\infty)$ with density $f^{\alpha,\lambda}$ defined in \eqref{eq:expressionoff}.}
  {Moreover,}
  \[
    {\sup_{t\in[0,1]}\left|F^T(t)-F^{\alpha,\lambda}(t)\right|\to 0,}
  \]
  {where $F^{\alpha,\lambda}$ is defined in \eqref{eq:Falphalambda}.}
\end{lemma}

\begin{proof}
  We provide the proof in Section~\ref{subsec:proofoflaplaceofrho}.
\end{proof}

To prove Lemma \ref{lemma:laplaceofrho}, we will use the following discrete-type extended Karamata--Tauberian theorem.

\begin{lemma}[Discrete-type extended Karamata-Tauberian Theorem]\label{lemma:discreteKaramataTauberian}
  Let $(\beta_n)_{n\ge0}$ be a non-negative sequence, $c \ge 0$, $\rho > -1$, and let $l$ be a slowly varying function. For $s>0$, define the transform
  \[
    \hat{\beta}(s) := \sum_{n=0}^\infty (e^{-sn} - e^{-s(n+1)}) \beta_n = (1 - e^{-s}) \sum_{n=0}^\infty e^{-sn} \beta_n.
  \]
  Assume this series converges for all $s > 0$. If
  \[
    \beta_n \sim \frac{c n^\rho l(n)}{\Gamma(1+\rho)} \quad \text{as } n \to \infty,
  \]
  then
  \[
    \hat{\beta}(s) \sim c s^{-\rho} l(1/s) \quad \text{as } s \downarrow 0.
  \]
\end{lemma}

\begin{proof}
  We provide the proof in Section \ref{subsec:proofofdiscreteKaramataTauberian}.
\end{proof}

We now record the tightness inputs.

\begin{lemma}\label{lemma:tightnessofyandlambda}
  The sequences $(Y^T)$ and $(\Lambda^T)$ are $C$-tight.
\end{lemma}

\begin{proof}
  We provide the proof in Section \ref{subsec:proofoftightnessofyandlambda}.
\end{proof}

{
In the tightness proof we will repeatedly bound expressions involving the resolvent coefficients,
in particular sums of the form $\sum_{k\le s}(A_k^T)^2$. While Lemma~\ref{lemma:laplaceofrho} provides a
Laplace/weak-convergence characterization of the rescaled resolvent measure (useful later for identifying the limit),
it does not imply pointwise bounds for $A_k^T$. We therefore establish below a deterministic polynomial upper bound,
uniform in $T$, which will be the only pointwise input needed in the tightness argument.

\begin{lemma}[Uniform polynomial bound for the resolvent]\label{lem:Ak_uniform_bound}
Let $(\eta_n)_{n\ge1}$ satisfy Assumption~\ref{assumption:1}(i) with $\alpha\in(1/2,1)$ and $\sum_{n\ge1}\eta_n=1$.
Define the (critical) renewal/resolvent sequence
\[
A_k:=\sum_{m=1}^{\infty}(\eta^{*m})_k,\qquad k\ge1.
\]
Then there exists a constant $C>0$ such that for all $k\ge1$,
\[
A_k \le C\,k^{\alpha-1}.
\]
Consequently, for every $T$ and every $k\ge1$,
\[
A_k^T=\sum_{m=1}^{\infty}(\eta^{T*m})_k=\sum_{m=1}^{\infty}a_T^m(\eta^{*m})_k\le A_k\le C\,k^{\alpha-1}.
\]
\end{lemma}

\begin{proof}
  We provide the proof in Section~\ref{pf:Ak_uniform_bound}.
\end{proof}
}

Since $(Y^T)$ and $(\Lambda^T)$ are pure jump processes, the classical Kolmogorov-Chentsov tightness criterion is not applicable. To establish the tightness of $(Z^{T},Y^{T})$, we instead apply a recent criterion from \cite{horst2023convergence}, which we restate as Lemma \ref{xu2023} below for the reader's convenience. This result gives a new criterion to verify the $C$--tightness for a sequence of c\`adl\`ag processes and will be used in proving Lemma \ref{lemma:tightnessofyandlambda}. Additionally, proving the convergence of the sharp bracket of $Z^T$ in the discrete-time case presents more challenges than in the continuous-time case discussed by \cite{jaisson2016rough}. For instance, in a Hawkes process $(N_t)_{t \ge 0}$ with {intensity function}
\[
  \lambda_t = \nu + \int_0^t h(t-s) \, dN_s,
\]
the angle bracket of $(N_t)_{t \ge 0}$ is given by $\int_0^\cdot \lambda_s \, ds$, and the sharp bracket is simply $N_\cdot$, which is straightforward to derive. However, achieving a similar result in the discrete-time setting is not as immediate and requires more nuanced handling. Indeed, the value $\lambda_n$ in \eqref{process_definition} is analogous to the intensity process $\lambda_t$ in continuous-time Hawkes processes. In the proof section, we work directly with $(Y^T)$, {and Section}
\ref{subsec:Heuristic} provides an intuitive interpretation of the rescaled $\lambda^T$. {To make the proof clear, we recall the following result.}

\begin{lemma}[Criterion for $C$-tightness, Lemma 3.5 in \cite{horst2023convergence}]\label{xu2023}
  Let $(X^{(n)})_{n\ge1}$ be a sequence of càdlàg stochastic processes, where each process $X^{(n)} = (X_t^{(n)})_{t \ge 0}$ is defined on a suitable probability space. Let $\Delta_h X_t^{(n)} := X_{t+h}^{(n)} - X_t^{(n)}$ denote the increment of the process $X^{(n)}$ over a time interval of length $h$. If for some $T>0$ the following three conditions hold:
  \begin{enumerate}
    \item 
    \[
      \sup_{n\ge1}\mathbb E\left[ \left|X^{(n)}_0\right|^\beta \right]<\infty
    \]
    for some $\beta>0$;
    \item \begin{align*}
      \sup_{k=0,1,\cdots,\lfloor Tn^\theta\rfloor}\sup_{h\in[0,1/n^\theta]}|\Delta_hX_{k/n^\theta}^{(n)}|\rightarrow 0\ \text{in probability as}\ n\rightarrow\infty,
    \end{align*}
    for some $\theta>2$;
    \item there exist some constants $C>0$, $p\ge1$, $m\in\{1,2,\cdots\}$, $\theta>2$ and pairs $\{(a_i,b_i)\}_{i=1,\cdots,m}$ satisfying
        \[
          a_i\ge0,\ b_i>0,\ \min_{1\le i\le m}\{b_i+a_i/\theta\}>1,
        \]
        such that for all $n\ge1$ and $h\in(0,1)$,
        \[
          \sup_{t\in[0,T]}\mathbb E\left[ \left|\Delta_hX_t^{(n)}\right|^p \right]\le C\cdot\sum_{i=1}^m\frac{h^{b_i}}{n^{a_i}}.
        \]
  \end{enumerate}
  Then the sequence $(X^{(n)})_{n\ge1}$ is $C$-tight on $[0, T]$.
\end{lemma}

\begin{lemma}\label{lemma:YLambdaconvergeinprobability}
  {
  The martingale part of the rescaled counting process vanishes: as $T\to\infty$,
  }
  \[
    {
    \sup_{t\in[0,1]}\left|Y_t^T-\Lambda_t^T\right|\xrightarrow{\mathbb P}0.}
  \]
  {
  Equivalently, $Y^T-\Lambda^T$ converges to $0$ in probability uniformly on $[0,1]$.
  }
\end{lemma}

\begin{proof}
  We provide the proof in Section \ref{subsec:proofofYLambdaconvergeinprobability}.
\end{proof}

We now state the tightness and quadratic variation result for the sequence $(Z^T,Y^T)_{T\ge 1}$.

\begin{proposition}[Tightness of the Sequence]\label{proposition:1}
  Under Assumption \ref{assumption:1} and \ref{assumption:2}, the sequence $(Z^T,Y^T)$ is tight. Furthermore, if $(Z,Y)$ is a limit point of $(Z^T,Y^T)$ as $T\rightarrow\infty$, then $Z$ is a continuous martingale and $[Z,Z]=Y$.
\end{proposition}
\begin{proof}
  We provide the proof in Section \ref{subsec:proofofproposition}.
\end{proof}

\paragraph{Renewal equation and characterization of the limit.}To characterize the limit $Y$, we will repeatedly use the explicit solution of discrete renewal equations.

\begin{lemma}[Solution of Discrete Renewal Equations]\label{lemma:discreterenewalequation}
  Given a non-negative sequence $(\eta_n)_{n\ge1} \in \ell^1$ and two non-negative sequences $(x_n)_{n\ge1}$, $(y_n)_{n\ge1}$, the following equation
  \begin{equation}\label{seriesconvolution}
    x_n = y_n + \sum_{s=1}^{n-1} \eta_s x_{n-s}
  \end{equation}
  has the unique solution
  \[
    x_n = \left(y + y * A\right)(n) = y_n + \sum_{i=1}^{n-1} A_i y_{n-i},
  \]
  {where $A_n=\sum_{k=1}^\infty (\eta^{*k})_n$.}
\end{lemma}

\begin{proof}
    We provide the proof in Section \ref{subsec:proofofdiscreterenewalequation}.
\end{proof}

{\begin{remark}[Mean intensity and the renewal sequence]\label{rem:mean_lambda}
Recall that $\lambda_s^T=\mu^T+\sum_{k=1}^{s-1}\eta_k^T X_{s-k}^T$. Taking expectations and using $\mathbb E[X_{s-k}^T]=\mathbb E[\lambda_{s-k}^T]$ yields the renewal-type equation
\[
q_s^T:=\mathbb E[\lambda_s^T]
=\mu^T+\sum_{k=1}^{s-1}\eta_k^T q_{s-k}^T,\qquad s\ge1.
\]
By Lemma~\ref{lemma:discreterenewalequation}, this admits the explicit solution
\[
q_s^T=\mu^T+\mu^T\sum_{u=1}^{s-1}A_u^T,\qquad s\ge1,
\]
where $(A_u^T)_{u\ge1}$ is defined in \eqref{ATexpression}.
\end{remark}}

\begin{theorem}[Characterization of the Limit Process]\label{theorem:1}
  For any limit point $(Z,Y)$ defined in Proposition \ref{proposition:1}, {there exists a Brownian motion $B$ (possibly on an extension of the underlying probability space) such that for all $t \in [0,1]$,}
  {$Z_t = B_{Y_t}$.}
  Furthermore, for any $\epsilon > 0$, $Y$ is continuous with H\"older regularity $(1 \wedge 2\alpha) - \epsilon$ on $[0,1]$ and satisfies
  \begin{equation}\label{eq:theorem1equation}
    Y_t = \int_0^t s f^{\alpha,\lambda}(t - s) \, ds + \frac{1}{\sqrt{\nu^* \lambda}} \int_0^t f^{\alpha,\lambda}(t - s) B_{Y_s} \, ds.
  \end{equation}
\end{theorem}

\begin{proof}
  The proof of the H\"older property for $Y$ in Theorem \ref{theorem:1} follows exactly the same argument as the Proof of Theorem 3.1 in \cite{jaisson2016rough}. The derivation of the equation \eqref{eq:theorem1equation} is provided in Section~\ref{subsec:proofoftheorem1}. 
\end{proof}

After \eqref{eq:theorem1equation} is established, we can characterize the derivative process (instantaneous variance in the financial interpretation).

\begin{corollary}[Characterization of the Derivative Process]\label{corollary:derivativesde}
  Let $(Y_t)$ be a process satisfying \eqref{eq:theorem1equation} for $t\in[0,1]$ and assume $\alpha>1/2$. Then $Y$ is differentiable on $[0,1]$ and its derivative $\dot{Y}$ satisfies
  \[
    \dot{Y}_t=F^{\alpha,\lambda}(t)+\frac1{\sqrt{\nu^*\lambda}}\int_0^tf^{\alpha,\lambda}(t-s)\sqrt{\dot{Y}_s}dB_s,
  \]
  where $B$ is the same Brownian motion as in \eqref{eq:theorem1equation}. Furthermore, for any $\epsilon>0$, $\dot{Y}$ has H\"older regularity $\alpha-1/2-\epsilon$.
\end{corollary}
\begin{proof}
The proof is exactly the same as the Proof of Theorem 3.2 in \cite{jaisson2016rough}. {More details on the proof are given in Section~\ref{subsec:proofofcorollaryinitialvalue}.}
\end{proof}

\subsection{Heuristic Motivation for the Scaling Limit}
\label{subsec:Heuristic}

In the ``light-tailed'' case, when the condition $\sum_{n=1}^\infty n\eta_n < \infty$ is met as in \cite{barczy2011asymptotic}, the scaling limit of the INAR process is known to be a squared Bessel process. Our work, however, focuses on the heavy-tailed case governed by Assumption \ref{assumption:1}, which induces a different scaling behavior.

To gain an intuitive understanding of this limit, we follow an approach inspired by the analysis of continuous-time Hawkes processes by \cite{jaisson2016rough}. The central idea is to study the asymptotic behavior of a properly renormalized version of the conditional intensity process, $(\lambda_n^T)$. Let us define the renormalized intensity process $(C_t^T)_{t\in[0,1]}$ as:
\begin{equation}\label{eq:renormalized_intensity}
  C_t^T := \frac{1-a_T}{\mu^T} \lambda_{\lfloor tT \rfloor}^T.
\end{equation}
The subsequent analysis in this section will heuristically show how this discrete process can be approximated by a continuous-time fractional Cox-Ingersoll-Ross process as $T \to \infty$. This provides an intuitive justification for the main results presented in our theorems.

To begin this heuristic analysis, we first introduce an alternative expression for $\lambda_n^T$ that will be useful in our derivations.
\begin{lemma}\label{lemma:anotherexpressionoflambda}
  For every $n\in \mathbb N$, 
  \begin{equation}\label{lambdaanotherexpression}
    \lambda_n^T=\mu^T+\mu^T\sum_{s=1}^{n-1}A_{n-s}^T+\sum_{s=1}^{n-1}A_{n-s}^T\left( M_{s}^T-M_{s-1}^T \right).
  \end{equation}
\end{lemma}

\begin{proof}
  We provide the proof in Section \ref{subsec:proofofanotherexpressionoflambda}.
\end{proof}

We substitute \eqref{lambdaanotherexpression} into \eqref{eq:renormalized_intensity} and obtain
\begin{equation}\label{Cexpression}
  \begin{aligned}
  C_t^T
  =&\frac{(1-a_T)}{\mu^T}\lambda_{\lfloor tT \rfloor}^T\\
  =& (1-a_T)+(1-a_T)\sum_{s=1}^{\lfloor tT \rfloor-1}A_s^T+\frac{(1-a_T)}{\mu^T}\sum_{s=1}^{\lfloor tT \rfloor-1}A_{\lfloor tT \rfloor-s}^T(M_s^T-M_{s-1}^T).
  \end{aligned}
\end{equation}

{
The argument below is intended as intuition only. In particular, Lemma~\ref{lemma:laplaceofrho} provides convergence of Laplace transforms (i.e.\ weak convergence at the level of measures), which does not imply convergence of densities or pointwise asymptotics for $A_k^T$.
Nevertheless, guided by the Laplace-transform asymptotics in Lemma~\ref{lemma:laplaceofrho}, one may \emph{formally} expect that for large $T$ and $k$ with $k/T$ of order one,
\[
  A_k^T \approx \frac{\|A^T\|_1 v_T}{T}\left(\frac{k}{T}\right)^{\alpha-1}E_{\alpha,\alpha}\!\left(-v_T\left(\frac{k}{T}\right)^\alpha\right).
\]
With this heuristic approximation in mind, the sum may be viewed as a Riemann sum:
\begin{align*}
  \sum_{s=1}^{\lfloor tT \rfloor-1}A_s^T
  &= \sum_{s=1}^{\lfloor tT\rfloor-1}\frac{1}{T}\cdot (T A_s^T)
  \approx \int_0^t T A_{\lfloor Tx\rfloor}^T\,dx
  \approx \int_0^t \|A^T\|_1 v_T x^{\alpha-1}E_{\alpha,\alpha}(-v_T x^\alpha)\,dx.
\end{align*}
A fully rigorous treatment used later in the paper relies only on uniform bounds and weak convergence of the associated discrete measures (rather than density convergence); see the proof of the moment bounds in Section~\ref{subsec:proofoftightnessofyandlambda}.
}

Define 
\[
  B_t^T:=\frac1{\sqrt{T}}\sum_{s=1}^{\lfloor tT \rfloor-1}\frac{M_s^T-M_{s-1}^T}{\sqrt{\lambda_s^T}}.
\]
By the martingale central limit theorem (see Section 3 in \cite{hall1980martingale}), $B_t^T$ converges to a standard Brownian motion. The last term in \eqref{Cexpression} can then be interpreted as a discrete convolution against a singular kernel, which motivates the limiting Volterra structure in Theorem~\ref{theorem:1}.

\paragraph{Microstructural foundation for the initial value.}
{
In practical applications, particularly in financial models such as the rough Heston model, it is essential to be able to prescribe a strictly positive initial instantaneous variance, i.e.\ $\dot Y_0=v_0$. 
In the classical Heston stochastic volatility model (see \cite{heston1993closed}), the asset price is driven by an instantaneous variance process $(V_t)_{t\ge0}$, and the integrated variance $Y_t:=\int_0^t V_s\,ds$ plays a central role in option pricing and likelihood-based inference. 
The rough Heston model (see \cite{el2019characteristic}) keeps the same interpretation---$V$ as instantaneous variance and $Y$ as integrated variance---but replaces the Markovian mean-reversion in $V$ by a fractional/Volterra-type kernel, yielding a non-Markovian variance process.
In our scaling limit, the process $Y$ corresponds to the integrated variance, so that $\dot Y$ represents the instantaneous variance. 
A subtle point is that if one starts the limiting Volterra dynamics from a time-independent initial condition, the resulting trajectories satisfy $\dot Y_0=0$, which would prevent us from matching a given market/empirical level $v_0>0$. 
For this reason, we adopt a time-dependent initialization (equivalently, we add a suitable deterministic term at the origin) so that the limiting model allows an arbitrary prescribed value $\dot Y_0=v_0$.
}

Following the methodology of \cite{el2019characteristic} and its discrete-time adaptation in \cite{wang2025rough}, we introduce a \textbf{time-dependent baseline intensity} $\hat{\mu}_T(n)$.

We redefine the intensity as
\begin{equation}\label{eq:lambdaT_time_dep}
  \lambda_n^T = \hat{\mu}_T(n) + \sum_{s=1}^{n-1} \eta_{n-s}^T X_s^T, \quad n \in \mathbb{N},
\end{equation}
where the time-dependent baseline intensity is given by:

{
\begin{definition}[Time-Dependent Baseline for Initial Value]\label{def:time_dependent_mu}
Fix $\xi \ge 0$. Define the time-dependent baseline intensity by
\[
  \hat{\mu}_T(n)
  := \mu^T + \xi \mu^T \left(
    \frac{1}{1-a_T}\left( 1 - \sum_{s=1}^{n-1}\eta_s^T \right)
    - \sum_{s=1}^{n-1}\eta_s^T
  \right), \qquad n\in\mathbb N,
\]
where $(\eta_n^T)_{n\ge 1}$ is a (possibly $T$-dependent) kernel sequence satisfying
$\sum_{n=1}^\infty \eta_n^T = a_T \in (0,1)$.

This baseline is strictly positive. Indeed, letting $S_{n-1}^T:=\sum_{s=1}^{n-1}\eta_s^T$, we have
$0\le S_{n-1}^T\le a_T$, hence
\[
  \frac{1}{1-a_T}(1-S_{n-1}^T)-S_{n-1}^T
  \ge \frac{1}{1-a_T}(1-a_T)-a_T
  = 1-a_T>0,
\]
and therefore
\[
  \hat{\mu}_T(n)\ge \mu^T+\xi\mu^T(1-a_T)>0,\qquad \forall\,n\ge1.
\]
\end{definition}

\noindent\textbf{Heuristic.}
At the very beginning ($n=1$), the sums are empty and
\[
  \hat{\mu}_T(1) = \mu^T\left(1+\frac{\xi}{1-a_T}\right),
\]
which is large when $(1-a_T)^{-1}$ is of order $\mathcal O(T^\alpha)$, creating an initial burst of events.
As $n$ increases, $S_{n-1}^T:=\sum_{s=1}^{n-1}\eta_s^T$ increases towards $a_T$, and the bracketed term
\[
  \frac{1}{1-a_T}(1-S_{n-1}^T)-S_{n-1}^T
\]
decreases from $\frac{1}{1-a_T}$ down to $1-a_T$.
Consequently, $\hat{\mu}_T(n)$ relaxes towards the positive asymptotic level
\[
  \lim_{n\to\infty}\hat{\mu}_T(n)
  = \mu^T+\xi\mu^T(1-a_T)
  = \mu^T\left(1+\xi(1-a_T)\right),
\]
rather than $\mu^T$.
In the scaling regime i.e. $(1-a_T)\to 0$, so this elevation is typically of lower order and
$\hat{\mu}_T(n)$ may be viewed as approaching $\mu^T$ up to a vanishing correction.
This engineered transient is what produces a non-zero initial value for the limiting \emph{rate} process $\dot Y_t$.

\medskip
With this modification, Theorem~\ref{theorem:1} extends to Corollary~\ref{cor:initial_value} below, showing that the
limit $\dot Y_t$ starts from a value determined by $\xi$.
We first introduce a lemma that will be used in the proof of Corollary~\ref{cor:initial_value}.}

\begin{lemma}[Renewal Expansion with Time-Dependent Baseline]\label{lemma:renewal_expansion_with_initial_value}
  Under the time-dependent baseline intensity $\hat{\mu}_T(n)$ from Definition \ref{def:time_dependent_mu}, the cumulative intensity and the centered cumulative process have the following representations for $n \ge 1$:
  \begin{equation}\label{eq:expectation_cumulative_lambda_initial_value}
    \mathbb{E}[N_n^T] = \sum_{i=1}^n \left( \hat{\mu}_T(i) + \sum_{s=1}^{i-1} A_{i-s}^T \hat{\mu}_T(s) \right),
  \end{equation}
  and
  \begin{equation}\label{eq:centered_N_initial_value}
    N_n^T - \mathbb{E}[N_n^T] = M_n^T + \sum_{s=1}^{n-1} A_{n-s}^T M_s^T.
  \end{equation}
\end{lemma}

\begin{proof}
    We provide the proof in Section \ref{subsection:proof_renewal_expansion_with_initial_value}.
\end{proof}

\begin{corollary}[Characterization of the Derivative Process with Initial Process Value]\label{cor:initial_value}
  Under the assumptions of Theorem \ref{theorem:1} and with the baseline intensity from Definition \ref{def:time_dependent_mu}, the limit process $(\dot{Y}_t)_{t\in[0,1]}$ is the unique solution to the stochastic Volterra equation
  \begin{equation}\label{eq:sve_with_initial}
    \dot{Y}_t=F^{\alpha,\lambda}(t)+\xi\left(1-F^{\alpha,\lambda}(t)\right)+\frac1{\sqrt{\nu^*\lambda}}\int_0^tf^{\alpha,\lambda}(t-s)\sqrt{\dot{Y}_s}dB_s,
  \end{equation}
  which implies the initial condition $\dot{Y}_0 = \xi$.
\end{corollary}

\begin{proof}
  We provide the proof in Section~\ref{subsec:proofofcorollaryinitialvalue}. 
\end{proof}
\section{Numerical Validation and Application}
\label{sec:numerical_validation}


The theoretical results of this paper show that the scaling limit of the heavy-tailed INAR($\infty$) process is an integrated fractional Cox-Ingersoll-Ross (CIR) process. This class of processes, governed by stochastic Volterra equations with singular kernels, has become central to modern financial mathematics for modeling so-called ``rough volatility'', which captures the observed long-memory and low smoothness properties of market volatility time series.

However, the direct numerical simulation of such continuous-time fractional processes presents significant challenges. Standard discretization methods, such as the Euler-Maruyama scheme, often suffer from slow convergence rates and require extremely small time steps to achieve reasonable accuracy, making them computationally expensive.

\subsection{Simulation via INAR Approximation}
\label{subsec:inar_simulation}

Our main theoretical result, encapsulated in Theorem \ref{theorem:1} and Corollary \ref{cor:initial_value}, provides a novel pathway for simulating the fractional CIR process. The core idea is that we can generate sample paths of the limiting continuous-time process by simulating its discrete-time antecedent---the heavy-tailed cumulative INAR($\infty$) process---for a large but finite $T_{\text{steps}}$. This approach leverages the computational simplicity of simulating discrete count processes, which only requires iterative sampling from a Poisson distribution.

A key feature of our method is the ability to specify a non-zero initial value, $v_0$, for the fractional CIR process ($\dot{Y}_t$). This is not achieved by an ad-hoc adjustment post-simulation, but is instead built directly into the microstructural dynamics of the INAR($\infty$) process. As outlined in Definition \ref{def:time_dependent_mu}, we employ a time-dependent baseline intensity, $\hat{\mu}_T(n)$, where the parameter $\xi$ is set to the desired initial value $v_0$. This creates an initial burst of activity that correctly sets the starting point of the process in the continuous-time limit.

The complete simulation procedure is detailed in Algorithm \ref{algorithm_inar_sim}.

\begin{algorithm}[h]
  \caption{Simulation of a fractional CIR process via the cumulative INAR($\infty$) approximation}
  \label{algorithm_inar_sim}
  \begin{algorithmic}[1]
  \STATE \textbf{Input:} Target initial value $v_0$, parameters $\alpha$, $\lambda$, $\nu^*$, and number of time steps $T_{\text{steps}}$.
  
  \STATE \textbf{Set} model parameter $\xi = v_0$.
  
  \STATE \textbf{Calculate} INAR($\infty$) model parameters based on the scaling limit theory:
  \STATE \quad $K = 1 / \zeta(1 + \alpha)$
  \STATE \quad $\delta = K \cdot \Gamma(1 - \alpha) / \alpha$
  \STATE \quad $\mu^T = \nu^* \cdot \delta^{-1} \cdot T_{\text{steps}}^{\alpha - 1}$
  \STATE \quad $a_T = 1 - (\lambda \cdot \delta) / T_{\text{steps}}^\alpha$
  
  \STATE \textbf{Pre-compute} the base kernel $\eta_n = K / n^{1+\alpha}$ and its cumulative sum for efficiency.
  
  \STATE \textbf{Initialize} count array $X^T$ of size $T_{\text{steps}}+1$ to zeros.
  
  \FOR{$n = 1$ \textbf{to} $T_{\text{steps}}$}
      \STATE Calculate the time-dependent baseline intensity $\hat{\mu}_T(n)$ using the formula from Definition \ref{def:time_dependent_mu}.
      \STATE Calculate the self-exciting (convolution) part: $C_n = \sum_{s=1}^{n-1} (a_T \eta_{n-s}) X_s^T$.
      \STATE Calculate the total intensity: $\lambda_n^T = \hat{\mu}_T(n) + C_n$.
      \STATE Draw the count for the current step: $X_n^T \sim \text{Pois}(\max(\lambda_n^T, 0))$.
  \ENDFOR
  
  \STATE \textbf{Post-processing to obtain the limit process}:
  \STATE Compute the cumulative process: $N^T_n = \sum_{s=1}^n X_s^T$.
  \STATE Apply the scaling factor: $Y_t^T = \frac{1 - a_T}{T_{\text{steps}}^\alpha \nu^* \delta^{-1}} N^T_{\lfloor t T_{\text{steps}} \rfloor}$.
  \STATE Compute the numerical derivative to get the final process: $V_t = \dot{Y}_t^T = \frac{\Delta Y^T}{\Delta t}$.
  
  \STATE \textbf{Output:} The simulated path $(V_t)_{t \in [0,1]}$.
  \end{algorithmic}
\end{algorithm}

\begin{figure}[H]
  \centering
  \includegraphics[width=0.8\textwidth]{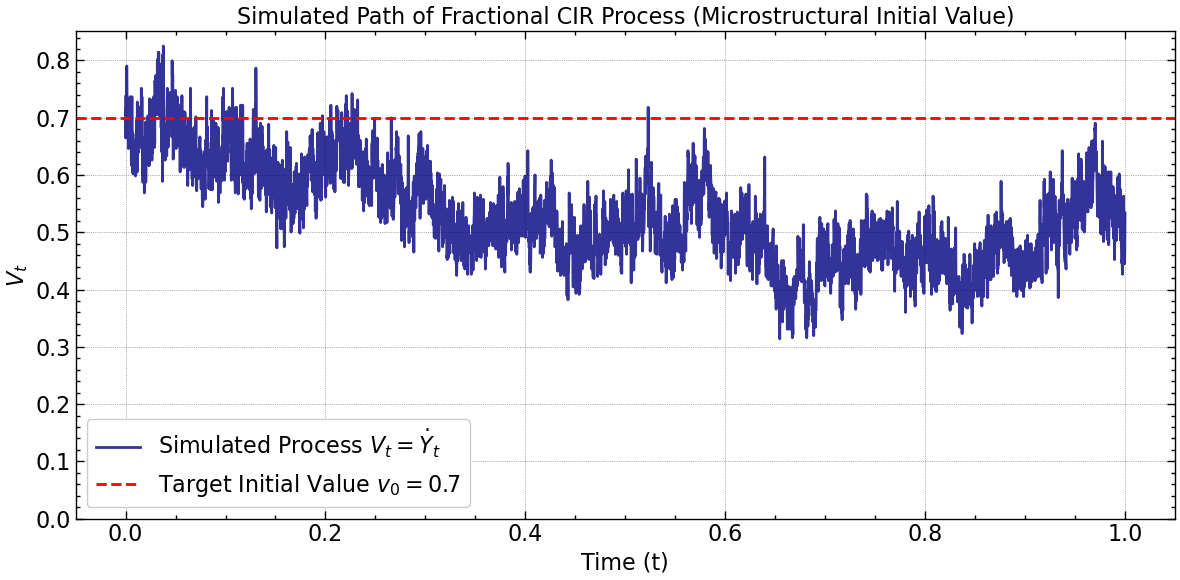}
  \caption{A sample path of the fractional CIR process $(V_t = \dot{Y}_t)_{0 \leq t \leq 1}$, simulated via the INAR($\infty$) approximation. The simulation parameters are $\alpha = 0.62$, $\lambda = 1.0$, $\nu^* = 200.0$, and a target initial value $v_0 = 0.7$. The number of discrete steps is $T_{\text{steps}}=5000$. The red dashed line indicates the target initial value, which the simulated path closely matches at $t=0$.}
  \label{fig:example_with_v0}
\end{figure}

The computational complexity of Algorithm \ref{algorithm_inar_sim} for generating a single path of length $T_{\text{steps}}$ is dominated by the main loop. Inside the loop, the calculation of the convolution term involves a sum of $n-1$ terms, resulting in a total complexity of $\mathcal{O}(T_{\text{steps}}^2)$ per simulation path. Figure \ref{fig:example_with_v0} displays a sample path generated using this algorithm, demonstrating its ability to start from the specified initial value.

\subsection{A Goodness-of-Fit Test and application in rough Heston models}
\label{subsec:goodness_of_fit}

Direct validation of the univariate fractional CIR process is challenging due to a scarcity of reliable benchmarks; therefore, we first highlight the framework's successful application to its well-benchmarked bivariate counterpart, the rough Heston model.
{
Recall that under the risk-neutral measure, the rough Heston model specifies the asset price $S$ and the instantaneous variance $V$ by
\[
  dS_t = S_t\sqrt{V_t}\,dW_t,
  \qquad
  V_t = V_0 + \frac{1}{\Gamma(\alpha)}\int_0^t (t-s)^{\alpha-1}\gamma\left(\theta - V_s\right)\,ds
        + \frac{1}{\Gamma(\alpha)}\int_0^t (t-s)^{\alpha-1}\gamma\nu\sqrt{V_s}\,dB_s,
\]
where $\alpha\in(1/2,1)$ controls the roughness (Volterra memory), $\gamma>0$ is the mean-reversion speed, $\theta>0$ the long-run variance level, $\nu>0$ the volatility-of-volatility, and $\rho\in(-1,1)$ is the correlation between the Brownian motions $(W,B)$.
For consistency with the benchmark parametrization used later in the numerical section, we employ the standard reparametrization in terms of $(\beta,\mu,\xi)$:
\[
\rho=\frac{1-\beta}{\sqrt{2(1+\beta^2)}}
\quad\Longleftrightarrow\quad
(2\rho^2-1)\beta^2+2\beta+(2\rho^2-1)=0,
\]
hence
\[
\beta=\frac{1\pm 2|\rho|\sqrt{1-\rho^2}}{1-2\rho^2}.
\]

\[
\mu=\frac{\theta(1+\beta^2)}{\gamma\nu^2(1+\beta)^2},
\qquad
\xi=\frac{V_0}{\theta}.
\]
so that the correlation/leverage and the variance-level inputs $(\rho,\theta,V_0)$ are matched exactly by the derived parameters $(\beta,\mu,\xi)$ reported below.}
To rigorously test the accuracy and goodness-of-fit of our INAR($\infty$) approximation, we focus on the rough volatility case by setting the parameter to $\alpha=0.62$. 

{We emphasize that our main theoretical results concern a univariate limit. However, for the univariate fractional CIR/Volterra variance dynamics there is, to the best of our knowledge, no closed-form or widely used high-accuracy benchmark that would allow a quantitative goodness-of-fit assessment comparable to the option-pricing benchmarks available for rough Heston. We therefore validate the INAR($\infty$) simulation methodology in the well-benchmarked rough Heston setting, where accurate reference prices are available (e.g.\ \cite{callegaro2021fast}). This numerical experiment is included as a consistency check for the microstructural INAR($\infty$) approximation and its simulation accuracy.}

{
\begin{itemize}
  \item \textit{Simulation setup.} We use the standard rough Heston parameter set from \cite{callegaro2021fast} (also used in \cite{wang2025rough}): $\gamma=0.1$, $\rho=-0.681$, $\nu=0.331$, $\theta=0.3156$, $V_0=0.0392$, $S_0=100$, and the corresponding derived parameters $\beta=27.5583$, $\mu=26.8592$, and $\xi=0.124208$. We simulate with $T=320$ time steps per year (i.e.\ $\Delta t = 1/T$).

  \item \textit{Bivariate microstructural setting (Table~\ref{tab:goodnessoffit}).} To simulate the rough Heston benchmark, we use the bivariate INAR($\infty$) microstructural model introduced in \cite{wang2025rough}. In that model, the reproduction coefficients are given by a matrix-valued kernel $(\eta_n^T)_{n\ge1}$ with $\eta_n^T\in\mathbb{R}_+^{2\times2}$, defined as $\eta_n^T=a_T\eta_n$. The baseline kernel $(\eta_n)_{n\ge1}$ is heavy-tailed (entrywise) and normalized so that the associated reproduction operator is critical (spectral radius equal to $1$), while $a_T\uparrow 1$ enforces near-instability as in Assumptions~\ref{assumption:1}--\ref{assumption:2}. We refer to \cite{wang2025rough} for the full construction linking the two-dimensional order-flow process to the rough Heston variance/price dynamics; here we only use this setting as a numerical validation benchmark for the INAR($\infty$) simulation methodology. This bivariate experiment does not affect the univariate limit theory developed earlier.

  \item \textit{Confidence intervals.} For each strike, the reported INAR($\infty$) price is the Monte Carlo mean of the (discounted) payoff over $M$ simulated paths, and the 95\% confidence interval is computed as
  \[
    \widehat{P}\ \pm\ 1.96\,\frac{\widehat{\sigma}}{\sqrt{M}},
  \]
  where $\widehat{\sigma}^2$ is the sample variance of the simulated discounted payoffs.
\end{itemize}
}

The robustness of the underlying INAR($\infty$) framework is highlighted by its successful application in related literature. In the work of \cite{wang2025rough}, for instance, a bivariate INAR($\infty$) process is developed to model the interacting dynamics of \textbf{buy and sell order flows}, providing a microstructural foundation for the full rough Heston model. The authors there prove the model's convergence and validate its accuracy through extensive numerical experiments. The success of the framework in that more complex, bivariate setting (as shown in \cite{wang2025rough}) provides strong collateral support for the validity of our approach. We therefore proceed with a direct numerical validation, comparing our simulation results against the benchmark from \cite{callegaro2021fast}.

\begin{table}[h]
\centering
\caption{Comparison of European option prices for the rough Heston model ($\alpha=0.62$) between the benchmark solution and our INAR($\infty$) simulation. The simulation uses 500,000 paths with {$T=320$}.}
\label{tab:goodnessoffit}
\begin{tabular}{ccc}
\hline
\textbf{Strike} & \textbf{Benchmark Price (\cite{callegaro2021fast})} & \textbf{INAR($\infty$) Simulation Price (95\% CI)} \\
\hline
80  & 22.1366 & 22.1558 [22.0972, 22.2145] \\
90  & 14.9672 & 14.9839 [14.9324, 15.0355] \\
100 & 9.4737  & 9.4898 [9.4470, 9.5326] \\
110 & 5.6234  & 5.6407 [5.6070, 5.6744] \\
120 & 3.1424  & 3.1584 [3.1331, 3.1837] \\
\hline
\end{tabular}
\end{table}

As shown in Table \ref{tab:goodnessoffit}, the prices generated by our INAR($\infty$) simulation for the rough case are in excellent agreement with the true prices from the benchmark solution. The differences are statistically insignificant and well within the 95\% confidence intervals. This provides strong numerical evidence for the accuracy and practical viability of our approximation method. The specific running code can be found in \url{https://github.com/gagawjbytw/INAR-rough-Heston}. 

While we use European options in this test for validation against a well-established benchmark, a key advantage of our simulation framework is its versatility. Being a Monte Carlo-based method, it naturally extends to the pricing of path-dependent options—such as Asian, Lookback, and Barrier options—for which analytical or semi-analytical solutions are often intractable. Moreover, the computational efficiency of this approach is a significant benefit. As demonstrated in the bivariate context by \cite{wang2025rough}, the INAR($\infty$) simulation framework is inherently parallelizable and proves to be considerably faster than traditional discretization schemes, such as the Euler method, applied to the limiting stochastic Volterra equation. It is stated in \cite{wang2025rough} that the total execution time on a {Mac mini M4 is 70.3283} seconds for calculating all four options (European, Asian, Lookback, Barrier) for {500,000} paths with {$T=320$}.


\section{Future Work}
In Section \ref{subsec:Heuristic}, we roughly used the scaling limit of the cumulative INAR($\infty$) process as an approximation of the integrated fractional CIR process. However, theoretically, the convergence rate of this scaling limit and the error analysis when using it as an approximation for the integrated fractional CIR process are still unknown. We will address this issue in future work. Another crucial direction for future research is the statistical estimation of the model's parameters. While the present paper focuses on the theoretical scaling limit, the practical application of our INAR($\infty$) framework requires robust methods for estimating the key parameters, such as those governing the heavy-tailed kernel, from observed data. This estimation problem is non-trivial due to the long-memory nature of the process. A comprehensive treatment of this topic, including the development and analysis of a Least-Square Estimation (LSE) procedure for this class of models, is provided by the authors in the separate work \cite{wang2024statistical}.


\section*{Acknowledgement}
The authors would like to thank the Editor, the Associate Editor, and two anonymous referees for their insightful and constructive feedback, which led to a substantial improvement in the quality and clarity of this paper. 
We are particularly indebted to the first referee for an exceptionally thorough and meticulous review; their detailed comments on our mathematical arguments were invaluable in helping us to significantly strengthen the rigor of our proofs. We are also deeply grateful to the second referee for their insightful, high-level comments, which prompted us to clarify the novelty and broader contribution of our work.
The authors are also grateful to Wei Xu, Yinhao Wu and Lingjiong Zhu for helpful comments.


\section*{Declarations}

\begin{itemize}
\item Conflict of interest/Competing interests: The authors declare that they have no known competing financial interests or personal relationships that could have appeared to influence the work reported in this paper.
\item Code availability: The code used in this study is provided as supplementary material. It is available under the \href{https://opensource.org/licenses/MIT}{MIT License} and can be freely accessed and used for non-commercial purposes.
\item Author contribution: Yingli Wang led the research, authored the main manuscript, provided the majority of the proofs, prepared Figure \ref{fig:example_with_v0}, and developed Algorithm \ref{algorithm_inar_sim}. Qinghua Wang contributed to the proof of Lemma \ref{lemma:tightnessofyandlambda} and played a key role in restructuring the manuscript. Ping He and Chunhao Cai supervised the project, provided critical revisions, and ensured the overall quality and coherence of the manuscript.
\end{itemize}


\section*{Funding Declarations}
Yingli Wang is supported by the Fundamental Research Funds for the Central Universities in Shanghai University of Finance and Economics CXJJ2023-397.

\bibliographystyle{chicago}
\bibliography{bibtex}

\appendix


\section{Proof of Lemmas}
\label{sec:proofoflemmas}


\subsection{Proof of Lemma \ref{lemma:sumandnotsum}}
\label{subsec:proofofsumandnotsum}

The proof consists of two main steps. First, we rigorously show that for a function $f(t) := t^{\beta} l(t)$, which is regularly varying with index ${\beta}$, its (partial) sum is asymptotically equivalent to the corresponding integral over unit intervals. Second, we apply Karamata's Integral Theorem(s).

{\medskip\noindent\textbf{Case (i): $\beta>-1$.}}
First, let's prove that $\sum_{t=X}^x f(t) \sim \int_X^x f(t)\,dt$. Our goal is to prove that
\[
\lim_{x\to\infty} \frac{\sum_{t=X}^x f(t)}{\int_X^x f(t)\,dt} = 1.
\]
We use the Stolz--Ces\`aro theorem. Let $a_x = \sum_{t=X}^x f(t)$ and $b_x = \int_X^x f(t)\,dt$. {Since $l$ is positive, $f(t)=t^\beta l(t)$ is positive for all $t\ge X$.} Hence $b_x$ is strictly increasing, and {since $\beta>-1$ we have $b_x\to\infty$ as $x\to\infty$}, satisfying the conditions for the theorem. We therefore analyze the limit of the ratio of differences:
\begin{equation*}
\lim_{x \to \infty} \frac{a_x - a_{x-1}}{b_x - b_{x-1}}
= \lim_{x \to \infty} \frac{f(x)}{\int_{x-1}^x f(t)\,dt}.
\end{equation*}
To evaluate this limit, we show that the denominator is asymptotically equivalent to the numerator. For any $t \in [x-1, x]$, we can write $t = \lambda x$, where $\lambda \in [1 - 1/x, 1]$. As $x \to \infty$, this interval for $\lambda$ converges to the point $\{1\}$.

According to \textbf{the Uniform Convergence Theorem for regularly varying functions} \cite[Theorem 1.5.2]{bingham1989regular}, the convergence of $f(\lambda x)/f(x) \to \lambda^{\beta}$ is uniform for $\lambda$ in any compact interval $[a,b]$ with $0<a \le b < \infty$. Let us choose such an interval that contains $1$, for example $[1/2, 2]$. For $x$ large enough, our interval $[1-1/x, 1]$ is contained within $[1/2, 2]$.

The uniform convergence implies that for any $\varepsilon > 0$, there exists an $N$ such that for all $x > N$ and for all $t \in [x-1, x]$ (and thus for all corresponding $\lambda \in [1-1/x, 1]$), we have:
\begin{equation*}
(1-\varepsilon)f(x) \le f(t) \le (1+\varepsilon)f(x).
\end{equation*}
Integrating this inequality over the interval $[x-1, x]$, which has length 1, yields:
\begin{equation*}
(1-\varepsilon)f(x) \le \int_{x-1}^x f(t)\,dt \le (1+\varepsilon)f(x).
\end{equation*}
This demonstrates that $\int_{x-1}^x f(t)\,dt \sim f(x)$ as $x \to \infty$. Therefore, the limit of the ratio of differences is 1. By the Stolz--Ces\`aro theorem, the limit of the original ratio is also 1, which proves that $\sum_{t=X}^x f(t) \sim \int_X^x f(t)\,dt$.

Next, we apply \textbf{Karamata's Integral Theorem}. Karamata's theorem for integrals of regularly varying functions \cite[Proposition 1.5.8]{bingham1989regular} states that for ${\beta}>-1$,
\[
    \int_X^x t^{\beta} l(t)\,dt \sim \frac{1}{{\beta}+1} x^{{\beta}+1} l(x) \quad \text{as } x \to \infty.
\]
Combining the above asymptotic equivalences, we obtain
\[
    \sum_{t=X}^x t^{\beta} l(t) \sim \frac{1}{{\beta}+1}x^{{\beta}+1}l(x),
\]
which proves item (i).

{\medskip\noindent\textbf{Case (ii): $\beta<-1$.}}
{In this case, $f(t)=t^\beta l(t)$ is eventually positive and integrable at infinity. We first show that}
\[
{
\sum_{t=x+1}^\infty f(t) \sim \int_x^\infty f(t)\,dt.}
\]
{Define $A_x:=\sum_{t=x+1}^\infty f(t)$ and $B_x:=\int_x^\infty f(t)\,dt$. Then $A_x\downarrow 0$ and $B_x\downarrow 0$ as $x\to\infty$. We again use the Stolz--Ces\'aro theorem applied to the decreasing sequences $(A_x)$ and $(B_x)$, and consider}
\begin{equation*}
{
\lim_{x\to\infty}\frac{A_x-A_{x+1}}{B_x-B_{x+1}}
= \lim_{x\to\infty}\frac{f(x+1)}{\int_x^{x+1} f(t)\,dt}.}
\end{equation*}
{As in Case (i), by writing $t=\lambda(x+1)$ for $t\in[x,x+1]$ and using the Uniform Convergence Theorem \cite[Theorem 1.5.2]{bingham1989regular}, we obtain $\int_x^{x+1} f(t)\,dt \sim f(x+1)$. Hence the above limit equals $1$, and Stolz--Ces\'aro yields $A_x/B_x\to 1$, i.e.\ $\sum_{t=x+1}^\infty f(t)\sim \int_x^\infty f(t)\,dt$.}

{Next, we apply Karamata's theorem for tails of integrals of regularly varying functions \cite[Proposition 1.5.10]{bingham1989regular}, which states that for $\beta<-1$,}
\[
{
\int_x^\infty t^\beta l(t)\,dt \sim -\frac{1}{\beta+1}\,x^{\beta+1}l(x)\quad\text{as }x\to\infty.}
\]
{Combining the last two asymptotic relations proves item (ii):}
\[
{
\sum_{t=x+1}^\infty t^\beta l(t) \sim -\frac{1}{\beta+1}\,x^{\beta+1}l(x).}
\]
This completes the proof.


\subsection{Proof of Lemma \ref{lemma:laplaceofrho}}
\label{subsec:proofoflaplaceofrho}
  {
  First, note that $\rho^T$ is a non-negative step function. Since the renewal (resolvent) sequence was
  originally defined only for indices $n\ge1$ by
  \[
  A_n^T:=\sum_{m=1}^\infty(\eta^{T*m})_n,\qquad n\ge1,
  \]
  we now extend it to $n=0$ by the convention
  \[
  \eta_0:=0,\qquad A_0^T:=0,
  \]
  {solely to write the step-function formula for $\rho^T$ uniformly over all intervals
  $[k/T,(k+1)/T)$ with $k\ge0$.} With this convention, $\rho^T$ is constant on $[k/T,(k+1)/T)$ with value
  $T A_k^T/\|A^T\|_1$, where $\|A^T\|_1:=\sum_{n\ge1}A_n^T$. Therefore,
  \begin{align*}
  \int_0^\infty \rho^T(x)\,dx
  &=\sum_{k=0}^\infty \int_{k/T}^{(k+1)/T} T\frac{A_k^T}{\|A^T\|_1}\,dx
  =\sum_{k=0}^\infty \frac{A_k^T}{\|A^T\|_1} \\
  &=\frac{A_0^T+\sum_{k\ge1}A_k^T}{\|A^T\|_1}
  =\frac{0+\|A^T\|_1}{\|A^T\|_1}=1,
  \end{align*}
  so $\rho^T$ is indeed a probability density on $[0,\infty)$.
  }

  {We set $\eta_0:=0$ and consequently $A_0^T:=0$ for all $T$.} Since $\|\eta\|_1=1$ {and $\eta_n\ge0$ for all $n\ge1$}, $\eta$ can be seen as a discrete probability measure. {Moreover, for each $T$ we have $\eta_n^T=a_T\eta_n\ge0$, so $(\eta_n^T)_{n\ge1}$ is also a non-negative sequence of reproduction coefficients.}
  Assume $(G_i)_{i\ge1}$ is a sequence of i.i.d.\ random variables with distribution $\mathbb P(G_i=n)=\eta_n$ for $n\ge1$. We first consider the rescaled random sum
  \[
    J^T=\frac1{T}\sum_{i=1}^{I^T}G_i,
  \]
  where $I^T$ is a geometric random variable with parameter $1-a_T$. {Set $S^T:=\sum_{i=1}^{I^T}G_i\in\mathbb{N}$. Since $(G_i)$ are i.i.d.\ with pmf $(\eta_n)$, for $n\ge1$ we have}
  \begin{equation*}
  {
  \mathbb{P}(S^T=n)
  =\sum_{k=1}^\infty (1-a_T)a_T^{k-1}\,\mathbb{P}\!\left(\sum_{i=1}^k G_i=n\right)
  =\sum_{k=1}^\infty (1-a_T)a_T^{k-1}\,(\eta^{*k})_n.}
  \end{equation*}
  {Recalling that $\eta^T=a_T\eta$, this can be rewritten as}
  \begin{equation*}
  {
  \mathbb{P}(S^T=n)=\frac{1-a_T}{a_T}\sum_{k=1}^\infty (\eta^{T* k})_n.}
  \end{equation*}
  {By the definition of the renewal sequence $(A_n^T)$ (cf.\ \eqref{ATexpression}), we have $A_n^T=\sum_{k=1}^\infty (\eta^{T* k})_n$, and since $\sum_{n\ge1}(\eta^{T* k})_n=\|\eta^T\|_1^k=a_T^k$, it follows that $\|A^T\|_1=\sum_{k=1}^\infty a_T^k=a_T/(1-a_T)$. Hence}
  \begin{equation*}
  {
  \mathbb{P}\!\left(J^T=\frac{n}{T}\right)=\mathbb{P}(S^T=n)=\frac{A_n^T}{\|A^T\|_1},\qquad n\ge1.}
  \end{equation*}

  {Note that $J^T$ is discrete, whereas $\rho^T$ is a step \emph{density}. To connect them, let $U$ be uniformly distributed on $(0,1)$, independent of everything else, and define the ``jittered'' variable $\widetilde J^T:=J^T+U/T$. Then $\widetilde J^T$ is absolutely continuous on $[0,\infty)$ and, by construction, it has density $\rho^T$ as defined in \eqref{eq:rho_T_def}.}

  {To identify the limit law, we study the Laplace transform of $J^T$ (and thus of $\widetilde J^T$). Denote $\hat\rho^T(z):=\mathbb{E}[e^{-zJ^T}]$, $z\ge0$. It satisfies}
  \begin{equation}\label{laplaceofjn}
  \begin{aligned}
    \hat\rho^T(z)
    =&\mathbb E[e^{-zJ^T}]
    =\sum_{k=1}^\infty (1-a_T)a_T^{k-1}\mathbb E\left[e^{-\frac{z}{T}\sum_{i=1}^kG_i}\right]\\
    =&\sum_{k=1}^\infty (1-a_T)a_T^{k-1}\left( \hat\gamma\left(\frac{z}{T}\right) \right)^k
    =\frac{(1-a_T)\hat\gamma(z/T)}{1-a_T\hat\gamma(z/T)},
  \end{aligned}
  \end{equation}
  where $\hat\gamma$ represents the Laplace transform of $G_1$, i.e.,
  \[
    \hat\gamma(z)=\mathbb E[e^{-zG_1}]=\sum_{n=1}^\infty \eta_ne^{-zn}.
  \]

{Moreover, since $\widetilde J^T=J^T+U/T$ with $U\sim\mathrm{Unif}(0,1)$ independent, the Laplace transform of the density $\rho^T$ satisfies}
\begin{equation*}
{
\int_0^\infty e^{-zx}\rho^T(x)\,dx
=\mathbb{E}[e^{-z\widetilde J^T}]
=\mathbb{E}[e^{-zJ^T}]\,\mathbb{E}[e^{-zU/T}]
=\hat\rho^T(z)\,\frac{1-e^{-z/T}}{z/T}.}
\end{equation*}
{Since $(1-e^{-z/T})/(z/T)\to1$ as $T\to\infty$, the limiting Laplace transform for $\rho^T$ is the same as that obtained for $\hat\rho^T(z)$.}

By Abel's lemma, we obtain
\begin{align*}
  {\sum_{n=1}^\infty \eta_ne^{-zn}}
  =&{\sum_{n=0}^\infty\left( e^{-zn}-e^{-z(n+1)} \right)\sum_{k=1}^n \eta_k}\\
  =&{\sum_{n=0}^\infty \left( e^{-zn}-e^{-z(n+1)} \right)\left( 1-\sum_{k=n+1}^\infty \eta_k \right)}\\
  =&{1-\sum_{n=0}^\infty\left( e^{-zn}-e^{-z(n+1)} \right)\sum_{k=n+1}^\infty\eta_k}\\
  =&{1-(1-e^{-z})\sum_{n=0}^\infty e^{-zn}\sum_{k=n+1}^\infty \eta_k.}
\end{align*}
{By Assumption~\ref{assumption:1}(i) and Lemma~\ref{lemma:sumandnotsum}, we have}
\[
  {\sum_{k=n+1}^\infty \eta_k \sim \frac{K}{\alpha}\,n^{-\alpha}\quad\text{as}\quad n\to\infty.}
\] 
From Lemma \ref{lemma:discreteKaramataTauberian}, we obtain
\[
  \sum_{n=0}^\infty\left( e^{-zn}-e^{-z(n+1)} \right)\sum_{k=n+1}^\infty\eta_k\sim \frac{K}{\alpha}\Gamma(1-\alpha)z^\alpha, \quad \text{as } z \downarrow 0.
\]
This yields the asymptotic expansion for $\hat\gamma(z)$ around the origin:
\[
  \hat\gamma(z) = 1-K\frac{\Gamma(1-\alpha)}{\alpha}z^\alpha+o(z^\alpha), \quad \text{as } z \downarrow 0.
\]
Let us define 
\begin{equation}\label{deltav}
  v_T=\delta^{-1}T^\alpha(1-a_T),
\end{equation}
{where we recall that $\delta=K\frac{\Gamma(1-\alpha)}{\alpha}$ (see Assumption~\ref{assumption:2}).}
Under Assumption \ref{assumption:2}, we have
\[
  {\lim_{T\to\infty} v_T = \lambda.}
\]
Substituting the expansion of $\hat\gamma(z/T)$ into \eqref{laplaceofjn} and taking $T\to\infty$, we find that for any fixed $z \ge 0$:
\[
  {\lim_{T\to\infty} \hat\rho^T(z) = \frac{\lambda}{\lambda+z^\alpha}.}
\]
{By \eqref{eq:laplace_mittagleffler}, the right-hand side coincides with the Laplace transform of the limiting density $f^{\alpha,\lambda}$ (defined in \eqref{eq:expressionoff}).} By L\'evy's continuity theorem for Laplace transforms, the convergence of the transforms implies the weak convergence of the corresponding probability measures.

{
Since $P^T \stackrel{\mathrm{w}}{\longrightarrow} P$ and $F^{\alpha,\lambda}$ is continuous on $[0,\infty)$, we now justify the uniform convergence of the distribution functions on $[0,1]$. 
Recall that $F^T(t):=\int_0^t \rho^T(x)\,dx$ is the distribution function of $P^T$, hence $F^T$ is {nondecreasing, right-continuous, and satisfies $\lim_{t\downarrow 0}F^T(t)=0$ and $\lim_{t\uparrow\infty}F^T(t)=1$}. 
Therefore, {as a sequence of distribution functions converging weakly to a distribution function with continuous limit $F^{\alpha,\lambda}$}, P\'olya's theorem (see e.g. \cite[Chapter~4.3, Exercise~4]{chung2001course}) yields
\[
\sup_{t\in[0,1]}\left|F^T(t)-F^{\alpha,\lambda}(t)\right|\longrightarrow 0,\qquad T\to\infty.
\]
which completes the proof.}


\subsection{Proof of Lemma \ref{lemma:discreteKaramataTauberian}}
\label{subsec:proofofdiscreteKaramataTauberian}
  Our goal is to connect the asymptotic behavior of the sequence $(\beta_n)$ to its transform $\hat{\beta}(s)$. We achieve this via the sequence of partial sums, $V_n = \sum_{k=0}^{n-1} \beta_k$ (with $V_0=0$).

  First, let us establish the exact relationship between the transform of the sequence, $\hat{\beta}(s)$, and the transform of its partial sums, $\hat{V}(s)$. Following the definitions, we have:
  \begin{align*}
    \hat{V}(s) &= \sum_{n=0}^{\infty}(e^{-sn}-e^{-s(n+1)})\sum_{m=0}^{n-1}\beta_{m} \\
    &= \sum_{m=0}^{\infty} \beta_m \sum_{n=m+1}^{\infty}(e^{-sn}-e^{-s(n+1)}) \\
    &= \sum_{m=0}^{\infty}\beta_{m}e^{-s(m+1)} \\
    &= \frac{e^{-s}}{1-e^{-s}}\hat{\beta}(s).
  \end{align*}
  As $s \downarrow 0$, we have $e^{-s} \to 1$ and $1-e^{-s} \sim s$. This implies the asymptotic relationship $\hat{V}(s) \sim s^{-1}\hat{\beta}(s)$.

  Now, to find the asymptotic behavior of $\hat{V}(s)$, we first need the asymptotics of $V_n$. We are given that $\beta_k \sim \frac{c k^\rho l(k)}{\Gamma(1+\rho)}$. This is precisely the form required by Lemma \ref{lemma:sumandnotsum}. Applying the lemma, we find the asymptotic behavior of the partial sums:
  \[
    V_n = \sum_{k=0}^{n-1} \beta_k \sim \frac{c}{\Gamma(1+\rho)} \frac{n^{\rho+1}l(n)}{\rho+1} = \frac{c n^{\rho+1}l(n)}{\Gamma(\rho+2)}.
  \]
  {With the behavior of $V_n$ established, we apply Karamata's Tauberian theorem for power series (e.g., \cite[Corollary~1.7.3]{bingham1989regular}). Since
  \[
    V_n \sim \frac{c}{\Gamma(\rho+2)}\,n^{\rho+1}l(n),
  \]
  the Tauberian theorem yields}
  \[
    {\hat{V}(s) \sim \frac{c}{\Gamma(\rho+2)}\,\Gamma(\rho+2)\, s^{-(\rho+1)} l(1/s)
    = c\, s^{-(\rho+1)} l(1/s) \quad \text{as } s \downarrow 0.}
  \]
  Finally, combining this with our derived relationship $\hat{\beta}(s) \sim s\hat{V}(s)$, we obtain the desired result:
  \[
    \hat{\beta}(s) \sim c s^{-\rho} l(1/s) \quad \text{as } s \downarrow 0.
  \]
  This completes the proof.


\subsection{Proof of Lemma \ref{lemma:tightnessofyandlambda}}
\label{subsec:proofoftightnessofyandlambda}

We use the criterion from Lemma \ref{xu2023} to verify the $C$-tightness for the sequence of c\`adl\`ag processes $(Y_t^T)_{t\in[0,1]}$.

{Before verifying the conditions, we derive moment bounds that hold uniformly in the time index $s\in\{1,\dots,\lfloor T\rfloor\}$. When convenient, we also consider the macroscopic scaling $s=\lfloor tT\rfloor$ with $t\in(0,1]$.}
{We decompose}
\[
{
\lambda_s^T = q_s^T+\kappa_s^T,\qquad q_s^T:=\mathbb E[\lambda_s^T],\qquad \kappa_s^T:=\lambda_s^T-q_s^T,}
\]
{so that $\mathbb E[\kappa_s^T]=0$.}

\begin{itemize}
    \item We first calculate \textbf{the first Moment of $\lambda_s^T$}. {By Remark~\ref{rem:mean_lambda}, for all $s\ge1$,
    \[
    \mathbb E[\lambda_s^T]=q_s^T=\mu^T+\mu^T\sum_{u=1}^{s-1}A_u^T.
    \]
    Since $\sum_{u\ge1}A_u^T=\|A^T\|_1=\frac{a_T}{1-a_T}$, it follows that
    \[
    q_s^T\le \mu^T\left(1+\sum_{u\ge1}A_u^T\right)=\frac{\mu^T}{1-a_T}.
    \]
    Consequently, using Assumption~\ref{assumption:2}, there exists a constant $\bar C_1>0$ such that
    \begin{equation}\label{eq:suplambdasbound}
    \sup_{1\le s\le \lfloor T\rfloor}\mathbb E[\lambda_s^T]\le \bar C_1\,T^{2\alpha-1},\qquad T\ge1.
    \end{equation}
    }

    \item Next, we calculate \textbf{the second Moment of $\lambda_s^T$}. The second moment is $\mathbb E[(\lambda_s^T)^2] = (q_s^T)^2 + \mathbb E[(\kappa_s^T)^2]$. We now rigorously derive the order of $\mathbb E[(\kappa_s^T)^2]$. {Recall that $\kappa_s^T=\lambda_s^T-q_s^T$ is the centered part of the intensity.
    Using the recursion for $\lambda_s^T$ and the fact that $q_s^T=\mathbb E[\lambda_s^T]$, we can write $\kappa_s^T$ as a moving-average of centered innovations:
    \[
      \kappa_s^T=\sum_{k=1}^{s-1}A_k^T\,\xi_{s-k}^T,
      \qquad
      \xi_j^T:=X_j^T-\lambda_j^T .
    \]
    Conditionally on $\mathcal F_{j-1}^T$, we have $\mathbb E [\xi_j^T\mid\mathcal F_{j-1}^T]=0$ and
    $\mathrm{Var}(\xi_j^T\mid\mathcal F_{j-1}^T)=\mathrm{Var}(X_j^T\mid\mathcal F_{j-1}^T)=\lambda_j^T$
    (since $X_j^T\mid\mathcal F_{j-1}^T\sim\mathrm{Pois}(\lambda_j^T)$).
    Moreover, $(\xi_j^T)$ is a martingale difference sequence with respect to $(\mathcal F_j^T)$.
    Hence for $k\neq \ell$ (say $k>\ell$), by the tower property,
    \[
      \mathbb E[\xi_{s-k}^T\,\xi_{s-\ell}^T]
      =\mathbb E\!\left[\xi_{s-\ell}^T\,\mathbb E[\xi_{s-k}^T\mid \mathcal F_{s-\ell}^T]\right]=0,
    \]
    since $\mathbb E[\xi_{s-k}^T\mid \mathcal F_{s-k-1}^T]=0$ and $\mathcal F_{s-\ell}^T\subseteq \mathcal F_{s-k-1}^T$ when $k>\ell$.
    Furthermore, using $\mathbb E[\xi_j^T\mid\mathcal F_{j-1}^T]=0$, we have
    \[
      \mathbb E[(\xi_j^T)^2\mid\mathcal F_{j-1}^T]=\mathrm{Var}(\xi_j^T\mid\mathcal F_{j-1}^T)
      =\mathrm{Var}(X_j^T\mid\mathcal F_{j-1}^T)=\lambda_j^T,
    \]
    and therefore $\mathbb E[(\xi_j^T)^2]=\mathbb E[\lambda_j^T]$.
    Consequently,
    \[
      \mathbb E[(\kappa_s^T)^2]
      =\sum_{k=1}^{s-1}(A_k^T)^2\,\mathbb E[(\xi_{s-k}^T)^2]
      =\sum_{k=1}^{s-1}(A_k^T)^2\,\mathbb E[\lambda_{s-k}^T].
    \]
    }

    Our goal is to prove that $\mathbb{E}[(\kappa_s^T)^2] \le C T^{4\alpha-2}$ for some constant $C>0$ as $T\to\infty$ and $s/T \to t \in (0,1]$. 
    
    {
    First, we calculate an upper bound for $\sum_{k=1}^{s-1} (A_k^T)^2$ that is uniform in $s\le T$.

    Define the (finite) measure on $(0,\infty)$
    \[
    \widetilde\mu_T(dx):=\frac{T}{\|A^T\|_1}\sum_{k\ge1}A_k^T\,\delta_{k/T}(dx).
    \]
    Lemma~\ref{lemma:laplaceofrho} identifies the limit of the Laplace transforms
    $\int_0^\infty e^{-zx}\widetilde\mu_T(dx)$, which implies $\widetilde\mu_T\Rightarrow \widetilde\mu$ in the sense of weak convergence of finite measures.
    In particular, for any bounded continuous test function $\varphi$,
    \[
    \int_0^\infty \varphi(x)\,\widetilde\mu_T(dx)\longrightarrow \int_0^\infty \varphi(x)\,\widetilde\mu(dx).
    \]

    {
    We now bound $\sum_{k=1}^{s-1}(A_k^T)^2$ uniformly for $s\le T$.
    By Lemma~\ref{lem:Ak_uniform_bound}, we have the deterministic pointwise bound
    $A_k^T\le C k^{\alpha-1}$ for all $k\ge1$ and all $T$.
    Hence, for $s\le T$,
    \begin{align*}
    \sum_{k=1}^{s-1}(A_k^T)^2
    &\le C^2\sum_{k=1}^{s-1}k^{2\alpha-2}
    \le C^2\sum_{k=1}^{T}k^{2\alpha-2}
    \le C^2\left(1+\int_1^{T}x^{2\alpha-2}\,dx\right)\\
    &= C^2\left(1+\frac{T^{2\alpha-1}-1}{2\alpha-1}\right)
    \le \bar C_2\,T^{2\alpha-1},
    \end{align*}
    for some constant $\bar C_2>0$ independent of $T$.
    }
    }

    Next, we establish a non-asymptotic upper bound for $\mathbb{E}[(\kappa_s^T)^2]$.

    First, from our preliminary calculations, for all $T \ge 1$ and $j \ge 1$, there exists a uniform constant $\bar{C}_1$ such that $\mathbb{E}[\lambda_j^T] \le \bar{C}_1 T^{2\alpha-1}$. 

    Furthermore, the analysis in the previous paragraph establishes that there exists a constant $\bar{C}_2$, independent of both $t$ and $T$, such that for any $t \in [0,1]$ and all $T \ge 1$:
    \[ \sum_{k=1}^{\lfloor tT \rfloor - 1} (A_k^T)^2 \le \bar{C}_2 T^{2\alpha-1}. \]
    {
    (In the previous paragraph we obtained $\sum_{k=1}^{\lfloor tT\rfloor-1}(A_k^T)^2 \le C_A^2\,T^{2\alpha-1}\int_{0}^{1}\frac{u^{2\alpha-2}}{(1+u^\alpha)^2}\,du$; hence we may take $\bar C_2:=C_A^2\int_{0}^{1}\frac{u^{2\alpha-2}}{(1+u^\alpha)^2}\,du$, which is independent of $t$ and $T$.)
    }

    Applying these uniform bounds for $s = \lfloor tT \rfloor$, we get:
    \begin{align*}
        \mathbb{E}[(\kappa_s^T)^2] 
        &= \sum_{k=1}^{s-1} (A_k^T)^2 \mathbb{E}[\lambda_{s-k}^T] \\
        &\le \left( \sup_{j \ge 1} \mathbb{E}[\lambda_j^T] \right) \left( \sum_{k=1}^{s-1} (A_k^T)^2 \right) \\
        &\le (\bar{C}_1 T^{2\alpha-1}) \cdot (\bar{C}_2 T^{2\alpha-1}) \\
        &= \bar{C}_1 \bar{C}_2 T^{4\alpha-2}.
    \end{align*}
    This provides the required non-asymptotic upper bound for the second moment.

    Since $(q_s^T)^2 =(\mathbb E[\lambda_s^T])^2 \le (\bar C_1)^2T^{4\alpha-2}$, we finally conclude that 
    \[
      \mathbb E[(\lambda_s^T)^2] \le \left((\bar C_1)^2+\bar C_1\bar C_2\right)\cdot T^{4\alpha-2}.
    \]
    
    \item \textbf{Second Moment of $X_s^T$}: From $\mathbb E[(X_s^T)^2] = \mathbb E[\lambda_s^T] + \mathbb E[(\lambda_s^T)^2]$, we have 
    \begin{equation}\label{eq:XTsecondmoment}
        \mathbb E[(X_s^T)^2] 
        \le \bar C_1 T^{2\alpha-1}+\left((\bar C_1)^2+\bar C_1\bar C_2\right) T^{4\alpha-2}
        \le \left(\bar C_1+(\bar C_1)^2+\bar C_1\bar C_2\right) T^{4\alpha-2}
    \end{equation}
    for $T\ge1$.
\end{itemize}

    {
    Now let us verify the three conditions in Lemma~\ref{xu2023}, which provide a sufficient criterion for $C$-tightness of the sequence $(Y^T)_{T\ge1}$.
    }
    \begin{enumerate}
        \item \textbf{Condition 1 (Initial Value):} Since $Y_0^T=0$, this condition is trivially satisfied.

        {
              \item \textbf{Condition 2 (Small increments in probability):}
              Fix $\theta>2$ and consider $t=k/T^\theta$ and $h\in[0,T^{-\theta}]$ as in Lemma~\ref{xu2023}.
              Since $Th\le T^{1-\theta}<1$ for all $T\ge2$, we have the deterministic floor-function bound
              \[
              {
              0\le \lfloor T(t+h)\rfloor-\lfloor Tt\rfloor \le 1.}
              \]
              {
              Hence the increment interval $(\lfloor Tt\rfloor,\lfloor T(t+h)\rfloor]$ contains either no integer or exactly one integer, and therefore}
              \[
              {
              \Delta_h Y_t^T
              = \frac{1-a_T}{T^\alpha \nu^*\delta^{-1}}
              \sum_{s=\lfloor Tt\rfloor+1}^{\lfloor T(t+h)\rfloor} X_s^T}
              \quad{\text{is either }0\text{ or equals } \frac{1-a_T}{T^\alpha \nu^*\delta^{-1}}\,X_{\lfloor T(t+h)\rfloor}^T.}
              \]
              {
              In particular, uniformly over $k\in\{0,1,\dots,\lfloor T^\theta\rfloor\}$ and $h\in[0,T^{-\theta}]$,}
              \[
              {
              \sup_{k,h}\left|\Delta_h Y_{k/T^\theta}^T\right|
              \le \frac{1-a_T}{T^\alpha \nu^*\delta^{-1}}\,
              \sup_{1\le s\le \lfloor T\rfloor} X_s^T.}
              \]
              }
              Let $C_T(\epsilon) = \frac{T^\alpha \nu^*\delta^{-1}}{1-a_T}\epsilon = \mathcal{O}(T^{2\alpha})$. There exists a constant $\bar C_3>0$ such that $\frac{T^\alpha \nu^*\delta^{-1}}{1-a_T} \ge \bar C_3T^{2\alpha}$ for $T\ge1$.
              By the union bound and Markov's inequality:
              \begin{align*}
                \mathbb P\left(\sup_{s \le \lfloor T \rfloor} X_s^T > C_T(\epsilon)\right) 
                \le &\sum_{s=0}^{\lfloor T \rfloor} \frac{\mathbb E[ (X_s^T)^2 ]}{C_T(\epsilon)^2} \\
                \le &T \cdot \frac{\left(\bar C_1+(\bar C_1)^2+\bar C_1\bar C_2\right) T^{4\alpha-2}}{\bar C_3^2T^{4\alpha}} \\
                = &\frac{\bar C_1+(\bar C_1)^2+\bar C_1\bar C_2}{\bar C_3^2}T^{-1} \to 0.
              \end{align*}
              Thus, Condition 2 is satisfied.

        \item \textbf{Condition 3 (Moment Bound on Increments):} First, there exists a constant $\bar C_5>0$ such that for any $T\ge1$ and $t,h\in[0,1]$:
              \[
                (\lfloor T(t+h) \rfloor - \lfloor Tt \rfloor + 1)^2\le \bar C_5 T^2h^2.
              \]
              Then by the Cauchy-Schwarz inequality:
              \begin{align*}
                \mathbb E\left[ \left| \Delta_h Y_t^T\right|^2 \right]
                \le &\left(\frac{1-a_T}{T^\alpha \nu^* \delta^{-1}}\right)^2 ( \lfloor T(t+h) \rfloor-\lfloor Tt \rfloor+1 ) \sum_{i=\lfloor Tt \rfloor}^{\lfloor T(t+h) \rfloor}\mathbb E\left[ (X_i^T)^2 \right]\\
                \le &\left(\frac{1-a_T}{T^\alpha \nu^* \delta^{-1}}\right)^2 ( \lfloor T(t+h) \rfloor-\lfloor Tt \rfloor+1 )^2\max_{i=\lfloor Tt \rfloor,\dots,\lfloor T(t+h) \rfloor}\mathbb E\left[ (X_i^T)^2 \right].
              \end{align*}
              From \eqref{eq:XTsecondmoment}, $\mathbb E[(X_i^T)^2] \le \left(\bar C_1+(\bar C_1)^2+\bar C_1\bar C_2\right) T^{4\alpha-2}$. Hence,
              \[
                \mathbb E\left[ \left| \Delta_h Y_t^T\right|^2 \right]
                \le \frac{\bar C_5(\bar C_1+(\bar C_1)^2+\bar C_1\bar C_2)}{\bar C_3^2}h^2.
              \]
              This demonstrates the condition holds with $p=2, b_1=2, a_1=0$. Thus, Condition 3 is satisfied.
      \end{enumerate}

      {
      \paragraph{C-tightness of \texorpdfstring{$(\Lambda^T)$}{LambdaT}.}
      It remains to show that $(\Lambda^T)$ is $C$-tight.
      Recall that by Lemma~\ref{lemma:YLambdaconvergeinprobability},
      \begin{equation}\label{eq:YT_minus_LambdaT_ucp}
        \sup_{t\in[0,1]}\left|Y_t^T-\Lambda_t^T\right|\xrightarrow[T\to\infty]{\mathbb P}0.
      \end{equation}
      Since $(Y^T)$ is $C$-tight, for any $\varepsilon>0$ there exists a compact set $K\subset C([0,1])$
      such that $\inf_T\mathbb P(Y^T\in K)\ge 1-\varepsilon$.
      For $\delta>0$, let
      \[
        K^\delta:=\left\{x\in D([0,1]):\ \inf_{y\in K}\|x-y\|_\infty\le \delta\right\},
      \]
      which is relatively compact in $D([0,1])$.
      Then, for all $T$,
      \[
        \mathbb P(\Lambda^T\in K^\delta)
        \ge \mathbb P\left(Y^T\in K,\ \|Y^T-\Lambda^T\|_\infty\le \delta\right)
        \ge 1-\varepsilon-\mathbb P\left(\|Y^T-\Lambda^T\|_\infty>\delta\right).
      \]
      Letting $T\to\infty$ and using \eqref{eq:YT_minus_LambdaT_ucp} yields
      $\liminf_{T\to\infty}\mathbb P(\Lambda^T\in K^\delta)\ge 1-\varepsilon$.
      This proves tightness of $(\Lambda^T)$.
      Moreover, any limit point of $(\Lambda^T)$ must coincide with the (continuous) limit of $(Y^T)$
      because of \eqref{eq:YT_minus_LambdaT_ucp}, hence $(\Lambda^T)$ is in fact $C$-tight.
      }


\subsection{Proof of Lemma~\ref{lem:Ak_uniform_bound}}\label{pf:Ak_uniform_bound}
\begin{proof}
We first note that for each $T$ and $m\ge1$,
\[
(\eta^{T*m})_k=(a_T\eta)^{*m}_k=a_T^m(\eta^{*m})_k,
\]
and therefore
\[
A_k^T=\sum_{m\ge1}(\eta^{T*m})_k=\sum_{m\ge1}a_T^m(\eta^{*m})_k
\le \sum_{m\ge1}(\eta^{*m})_k = A_k,
\]
since $a_T\in(0,1)$ and all terms are nonnegative. Hence it suffices to prove that
$A_k\le C k^{\alpha-1}$.

To this end, define the Laplace (power-series) transforms
\[
\hat\eta(s):=\sum_{n\ge1}\eta_n e^{-sn},\qquad s>0,
\qquad\text{and}\qquad
\hat A(s):=\sum_{k\ge1}A_k e^{-sk},\qquad s>0.
\]
Because $\eta$ is a probability mass function, $\hat\eta(s)\in(0,1)$ for every $s>0$.
Using the convolution property and Fubini's theorem (all terms are nonnegative), we obtain
\begin{align*}
\hat A(s)
&=\sum_{k\ge1}\sum_{m\ge1}(\eta^{*m})_k e^{-sk}
=\sum_{m\ge1}\sum_{k\ge1}(\eta^{*m})_k e^{-sk}
=\sum_{m\ge1}\left(\sum_{n\ge1}\eta_n e^{-sn}\right)^m \\
&=\sum_{m\ge1}\hat\eta(s)^m
=\frac{\hat\eta(s)}{1-\hat\eta(s)}.
\end{align*}

Next we analyze $1-\hat\eta(s)$ as $s\downarrow0$. By the same Abel-summation argument as in
the proof of Lemma~\ref{lemma:laplaceofrho},
\[
\hat\eta(s)
=1-(1-e^{-s})\sum_{n=0}^\infty e^{-sn}\sum_{k=n+1}^\infty\eta_k.
\]
Under Assumption~\ref{assumption:1}(i), the tail satisfies
\[
\sum_{k=n+1}^\infty\eta_k \sim \frac{K}{\alpha}\,n^{-\alpha}\qquad(n\to\infty).
\]
Applying Lemma~\ref{lemma:discreteKaramataTauberian} with
\[
\beta_n:=\sum_{k=n+1}^\infty\eta_k,\quad
\rho:=-\alpha,\quad
c:=\frac{K}{\alpha}\Gamma(1-\alpha),
\]
(which is legitimate since $\alpha\in(\frac12,1)$ implies $\rho=-\alpha>-1$) yields
\[
(1-e^{-s})\sum_{n=0}^\infty e^{-sn}\sum_{k=n+1}^\infty\eta_k
\sim \frac{K}{\alpha}\Gamma(1-\alpha)\,s^\alpha
=\delta\,s^\alpha
\qquad(s\downarrow0),
\]
where $\delta:=K\Gamma(1-\alpha)/\alpha$. Consequently,
\[
1-\hat\eta(s)\sim \delta\,s^\alpha\qquad(s\downarrow0).
\]
Combining this with $\hat\eta(s)\to1$ as $s\downarrow0$, we deduce
\[
\hat A(s)=\frac{\hat\eta(s)}{1-\hat\eta(s)}
\sim \frac{1}{\delta}\,s^{-\alpha}\qquad(s\downarrow0).
\]

Now define the partial sums
\[
U_n:=\sum_{k=1}^n A_k,\qquad n\ge1.
\]
A summation-by-parts identity gives
\[
(1-e^{-s})\sum_{n=0}^\infty e^{-sn}U_n
=\sum_{n=0}^\infty (e^{-sn}-e^{-s(n+1)})U_n
=\sum_{k\ge1}A_k e^{-sk}
=\hat A(s),
\]
and hence
\[
(1-e^{-s})\sum_{n=0}^\infty e^{-sn}U_n
\sim \frac{1}{\delta}\,s^{-\alpha}\qquad(s\downarrow0).
\]
Since $(U_n)$ is nondecreasing, we may apply the discrete Karamata Tauberian theorem for power series
(e.g.\ \cite[Corollary~1.7.3]{bingham1989regular}) to obtain
\[
U_n \sim \frac{1}{\delta\,\Gamma(1+\alpha)}\,n^\alpha\qquad(n\to\infty).
\]
In particular, there exist $C_0>0$ and $n_0$ such that $U_n\le C_0 n^\alpha$ for all $n\ge n_0$.
Therefore, for $n\ge n_0$,
\[
A_n = U_n-U_{n-1}\le C_0\left(n^\alpha-(n-1)^\alpha\right)
\le C_0\,\alpha\,n^{\alpha-1},
\]
where we used the mean-value theorem (equivalently, the concavity of $x\mapsto x^\alpha$ on $(0,\infty)$).
Enlarging the constant to also cover $1\le n<n_0$ yields
\[
A_n\le C\,n^{\alpha-1}\qquad\text{for all }n\ge1.
\]

Finally, since $A_k^T\le A_k$ for all $k$, we conclude that
\[
A_k^T\le A_k\le C\,k^{\alpha-1},\qquad \forall\,T,\ \forall\,k\ge1,
\]
which completes the proof.
\end{proof}



\subsection{Proof of Lemma \ref{lemma:YLambdaconvergeinprobability}}
\label{subsec:proofofYLambdaconvergeinprobability}
  We aim to show that for any $\epsilon > 0$, $\mathbb{P}(\sup_{t\in[0,1]} |Y_t^T - \Lambda_t^T| > \epsilon) \to 0$ as $T \to \infty$.
  We have
  \[
    Y_t^T-\Lambda_t^T = \frac{1-a_T}{T^\alpha \nu^*\delta^{-1}} M_{\lfloor tT \rfloor}^T,
  \]
  where $M_n^T = \sum_{s=1}^{n}(X_s^T-\lambda_s^T)$ is a zero-mean martingale.
  
  {
  Applying Doob's $L^2$-inequality to the martingale $(Y_t^T-\Lambda_t^T)_{t\in[0,1]}$ (equivalently to $M^T$), we get:
  }
  \begin{align*}
    \mathbb E\left[ \sup_{t\in[0,1]}\left( Y_t^T-\Lambda_t^T \right)^2 \right]
    \le 4 \left( \frac{1-a_T}{T^\alpha \nu^*\delta^{-1}} \right)^2 \mathbb E\left[(M_{\lfloor T\rfloor}^T)^2\right]\
  \end{align*}
  The second moment of the martingale $M_{\lfloor T\rfloor}^T$ is equal to the expectation of its predictable quadratic variation:
  \begin{align*}
    \mathbb E\left[(M_{\lfloor T\rfloor}^T)^2\right] &= \mathbb E\left[ \left( \sum_{s=1}^{\lfloor T\rfloor}(X_s^T-\lambda_s^T) \right)^2 \right] \\
    &= \mathbb E\left[\sum_{s=1}^{\lfloor T\rfloor}(X_s^T-\lambda_s^T)^2\right] && \text{(Cross-terms vanish for martingale differences)} \\
    &= \sum_{s=1}^{\lfloor T\rfloor} \mathbb E\left[ \mathbb E[(X_s^T-\lambda_s^T)^2 | \mathcal{F}_{s-1}] \right] \\
    &= \sum_{s=1}^{\lfloor T\rfloor} \mathbb E[\lambda_s^T] && \text{(Since } X_s^T|\mathcal{F}_{s-1} \sim \text{Pois}(\lambda_s^T)\text{)} \\
    &= \mathbb E\left[ \sum_{s=1}^{\lfloor T\rfloor} X_s^T \right] = \mathbb E[N_{\lfloor T\rfloor}^T]. && \text{(Since } \mathbb E[X_s^T] = \mathbb E[\lambda_s^T]\text{)}
  \end{align*}
  {We now bound $\mathbb E[N_{\lfloor T\rfloor}^T]$.
  Using the tower property and $\mathbb E[X_s^T\mid\mathcal F_{s-1}]=\lambda_s^T$, we have
  \[
  \mathbb E[N_{\lfloor T\rfloor}^T]
  =\sum_{s=1}^{\lfloor T\rfloor}\mathbb E[X_s^T]
  =\sum_{s=1}^{\lfloor T\rfloor}\mathbb E[\lambda_s^T]
  \le \lfloor T\rfloor \sup_{1\le s\le \lfloor T\rfloor}\mathbb E[\lambda_s^T].
  \]
  By the uniform first-moment bound \eqref{eq:suplambdasbound}, there exists $\bar C_1>0$ such that
  $\sup_{1\le s\le \lfloor T\rfloor}\mathbb E[\lambda_s^T]\le \bar C_1 T^{2\alpha-1}$ for all $T\ge1$.
  Therefore,
  \[
  \mathbb E[N_{\lfloor T\rfloor}^T]\le \bar C_1\,T\cdot T^{2\alpha-1}=\bar C_1 T^{2\alpha}.
  \]
  }

 {
  Substituting this bound back, we get:
  }
  \begin{align*}
  \mathbb E\left[ \sup_{t\in[0,1]}\left( Y_t^T-\Lambda_t^T \right)^2 \right]
  &\le 4 \left( \frac{1-a_T}{T^\alpha \nu^*\delta^{-1}} \right)^2 \mathbb E\left[(M_{\lfloor T\rfloor}^T)^2\right] \\
  &= 4 \left( \frac{1-a_T}{T^\alpha \nu^*\delta^{-1}} \right)^2 \mathbb E\left[N_{\lfloor T\rfloor}^T\right] \\
  &\le 4 \left( \frac{1-a_T}{T^\alpha \nu^*\delta^{-1}} \right)^2 \cdot \bar C_1 T^{2\alpha}
  \;\le\; C' (1-a_T)^2,
  \end{align*}
  {
  for some constant $C'>0$. Since $a_T\to 1$ by Assumption~\ref{assumption:2}, we have $(1-a_T)^2\to 0$, hence the right-hand side converges to $0$.
  }


\subsection{Proof of Lemma \ref{lemma:discreterenewalequation}}
\label{subsec:proofofdiscreterenewalequation}
  We consider the generating functions (or $z$-transforms) of the sequences involved. Since $(\eta_n)_{n\ge1}$ is a non-negative sequence in $\ell^1$, its generating function $\mathcal{G}(\eta)(z) = \sum_{n=1}^\infty \eta_n z^n$ converges absolutely for all complex $z$ with $|z| \le 1$. All subsequent manipulations are valid for $|z|<1$, where $|\mathcal{G}(\eta)(z)| < \sum_{n=1}^\infty \eta_n \le 1$.

  Taking the $z$-transform on both sides of equation \eqref{seriesconvolution} and using the convolution property \cite{oppenheim1999discrete}, we obtain:
  \[
    \mathcal{G}(x)(z) = \mathcal{G}(y)(z) + \mathcal{G}(\eta)(z) \cdot \mathcal{G}(x)(z).
  \]
  Solving for $\mathcal{G}(x)(z)$, we get:
  \[
    \mathcal{G}(x)(z) = \frac{\mathcal{G}(y)(z)}{1 - \mathcal{G}(\eta)(z)} = \mathcal{G}(y)(z) \cdot \left(1 + \frac{\mathcal{G}(\eta)(z)}{1 - \mathcal{G}(\eta)(z)}\right).
  \]
  To interpret the second term, we find the generating function of the sequence $(A_n)_{n\ge1}$. By definition, $A_n = \sum_{k=1}^\infty (\eta^{*k})_n$, where $\eta^{*k}$ is the $k$-fold convolution of $\eta$. The generating function of a convolution is the product of the generating functions. Thus, the generating function of $\eta^{*k}$ is $(\mathcal{G}(\eta)(z))^k$. The generating function of $A_n$ is therefore:
  \begin{align*}
    \mathcal{G}(A)(z) &= \sum_{n=1}^\infty A_n z^n = \sum_{n=1}^\infty \left( \sum_{k=1}^\infty (\eta^{*k})_n \right) z^n \\
    &= \sum_{k=1}^\infty \sum_{n=1}^\infty (\eta^{*k})_n z^n = \sum_{k=1}^\infty (\mathcal{G}(\eta)(z))^k.
  \end{align*}
  This is a geometric series, which for $|\mathcal{G}(\eta)(z)|<1$ converges to:
  \[
    \mathcal{G}(A)(z) = \frac{\mathcal{G}(\eta)(z)}{1 - \mathcal{G}(\eta)(z)}.
  \]
  Substituting this back into the expression for $\mathcal{G}(x)(z)$, we have:
  \[
    \mathcal{G}(x)(z) = \mathcal{G}(y)(z) \cdot (1 + \mathcal{G}(A)(z)) = \mathcal{G}(y)(z) + \mathcal{G}(y)(z)\mathcal{G}(A)(z).
  \]
  By performing the inverse $z$-transform and applying the convolution theorem again, we obtain the unique solution:
  \[
    x_n = y_n + \sum_{i=1}^{n-1} A_i y_{n-i}.
  \]
  This concludes the proof.

\subsection{Proof of Lemma \ref{lemma:anotherexpressionoflambda}}
\label{subsec:proofofanotherexpressionoflambda}
  We notice that
  \begin{align*}
    \lambda_n^T
    =&\mu^T+\sum_{s=1}^{n-1}\eta_{n-s}^TX_s^T
    =\mu^T+\sum_{s=1}^{n-1}\eta_{n-s}^T\left( M_s^T-M_{s-1}^T+\lambda_s^T \right)\\
    =&\mu^T+\sum_{s=1}^{n-1}\eta_{n-s}^T\left( M_s^T-M_{s-1}^T \right)+\sum_{s=1}^{n-1}\eta_{n-s}^T\lambda_s^T.
  \end{align*}
  By the Lemma \ref{lemma:discreterenewalequation}, we obtain
  \begin{align*}
    \lambda_n^T
    =&\mu^T+\sum_{s=1}^{n-1}\eta_{n-s}^T\left( M_s^T-M_{s-1}^T \right)+\sum_{s=1}^{n-1}A_{n-s}^T\left( \mu^T+\sum_{u=1}^s\eta_{s-u}^T\left( M_u^T-M_{u-1}^T \right)\right)\\
    =&\mu^T+\sum_{s=1}^{n-1}\eta_{n-s}^T\left( M_s^T-M_{s-1}^T \right)+\mu^T\sum_{s=1}^{n-1}A_{n-s}^T+\sum_{s=1}^{n-1}A_{n-s}^T\left( \sum_{u=1}^{s-1}\eta_{s-u}^T\left( M_u^T-M_{u-1}^T \right) \right)\\
    =&\mu^T+\sum_{s=1}^{n-1}\eta_{n-s}^T\left( M_s^T-M_{s-1}^T \right)+\mu^T\sum_{s=1}^{n-1}A_{n-s}^T+\sum_{u=1}^{n-1}(M_u^T-M_{u-1}^T)\sum_{s=u+1}^{n-1}\eta_{s-u}^TA_{n-s}^T\\
    =&\mu^T+\sum_{s=1}^{n-1}\eta_{n-s}^T\left( M_s^T-M_{s-1}^T \right)+\mu^T\sum_{s=1}^{n-1}A_{n-s}^T+\sum_{u=1}^{n-1}(M_u^T-M_{u-1}^T)\left( A_{n-u}^T-\eta_{n-u}^T \right)\\
    =&\mu^T+\mu^T\sum_{s=1}^{n-1}A_{n-s}^T+\sum_{s=1}^{n-1}A_{n-s}^T\left( M_{s}^T-M_{s-1}^T \right).
  \end{align*}

\subsection{Proof of Lemma \ref{lemma:renewal_expansion_with_initial_value}}
\label{subsection:proof_renewal_expansion_with_initial_value}

The proof relies on the application of the discrete renewal equation from Lemma \ref{lemma:discreterenewalequation}.

First, let's prove Equation \eqref{eq:expectation_cumulative_lambda_initial_value}. By taking the expectation and using the law of total expectation:
\begin{align*}
    \mathbb{E}[N_n^T] &= \sum_{i=1}^n \mathbb{E}[X_i^T] = \sum_{i=1}^n \mathbb{E}[\lambda_i^T] \\
    &= \sum_{i=1}^n \mathbb{E}\left[ \hat{\mu}_T(i) + \sum_{s=1}^{i-1} \eta_{s}^T X_{i-s}^T \right] \\
    &= \sum_{i=1}^n \hat{\mu}_T(i) + \sum_{i=1}^n \sum_{s=1}^{i-1} \eta_s^T \mathbb{E}[X_{i-s}^T] \\
    &= \sum_{i=1}^n \hat{\mu}_T(i) + \sum_{s=1}^{n-1} \eta_s^T \sum_{i=s+1}^n \mathbb{E}[X_{i-s}^T] \\
    &= \sum_{i=1}^n \hat{\mu}_T(i) + \sum_{s=1}^{n-1} \eta_s^T \mathbb{E}[N_{n-s}^T].
\end{align*}
This is a discrete renewal equation for $x_n = \mathbb{E}[N_n^T]$ with driving term $y_n = \sum_{i=1}^n \hat{\mu}_T(i)$. Applying the solution from Lemma \ref{lemma:discreterenewalequation}, we have:
\[
    \mathbb{E}[N_n^T] = \sum_{i=1}^n \hat{\mu}_T(i) + \sum_{s=1}^{n-1} A_s^T \left( \sum_{i=1}^{n-s} \hat{\mu}_T(i) \right).
\]
Re-arranging the order of summation in the second term gives the desired result.

Next, for Equation \eqref{eq:centered_N_initial_value}, we start by expressing $N_n^T$ in terms of its martingale and intensity components:
\[
    N_n^T = \sum_{i=1}^n X_i^T = \sum_{i=1}^n (M_i^T - M_{i-1}^T + \lambda_i^T) = M_n^T + \sum_{i=1}^n \lambda_i^T.
\]
Subtracting the expectation, we get:
\[
    N_n^T - \mathbb{E}[N_n^T] = M_n^T + \sum_{i=1}^n (\lambda_i^T - \mathbb{E}[\lambda_i^T]).
\]
Let $Z_n = N_n^T - \mathbb{E}[N_n^T]$. Then $\lambda_i^T - \mathbb{E}[\lambda_i^T] = \sum_{s=1}^{i-1} \eta_s^T (X_{i-s}^T - \mathbb{E}[X_{i-s}^T]) = \sum_{s=1}^{i-1} \eta_s^T Z_{i-s}$.
The equation becomes a discrete renewal equation for $Z_n$:
\[
    Z_n = M_n^T + \sum_{i=1}^n \sum_{s=1}^{i-1} \eta_s^T Z_{i-s} = M_n^T + \sum_{s=1}^{n-1} \eta_s^T Z_{n-s}.
\]
Applying the solution from Lemma \ref{lemma:discreterenewalequation} with driving term $y_n = M_n^T$ gives the result:
\[
    N_n^T - \mathbb{E}[N_n^T] = M_n^T + \sum_{s=1}^{n-1} A_{s}^T M_{n-s}^T = M_n^T + \sum_{s=1}^{n-1} A_{n-s}^T M_s^T.
\]
This completes the proof.


\section{Proof of the Main Results}
\label{sec:proofofthemainresults}

\subsection{Proof of Proposition \ref{proposition:1}}
\label{subsec:proofofproposition}

The proof proceeds in three main steps. First, we establish the joint tightness of the sequence $(Z^T, Y^T)$. Second, we show that the quadratic variation $[Z^T, Z^T]_t$ converges to the same limit as $Y^T_t$. Finally, we use this to characterize the limit process $(Z,Y)$.

\textbf{1. Tightness of $(Z^T, Y^T)$.}
From Lemma \ref{lemma:tightnessofyandlambda}, we have already established that the sequences $(Y^T)$ and $(\Lambda^T)$ are $C$-tight. The process $Z^T$ is defined as
\[
  Z_t^T = \sqrt{\frac{T^\alpha \nu^* \delta^{-1}}{1-a_T}} (Y_t^T - \Lambda_t^T).
\]
{
Since $(Y^T)$ and $(\Lambda^T)$ are $C$-tight in $\mathbb D([0,1])$, the pair $(Y^T,\Lambda^T)$ is $C$-tight in the product space $\mathbb D([0,1])^2$.
Moreover, $\mathbb D([0,1])$ endowed with the Skorokhod $J_1$ topology is a topological vector space, hence the subtraction map
\[
\Phi:\mathbb D([0,1])^2\to \mathbb D([0,1]),\qquad \Phi(x,y)=x-y,
\]
is continuous (see, e.g., \cite[Chapter~VI, \S1]{jacod2013limit}).
Therefore, by the continuous mapping theorem, $(Y^T-\Lambda^T)$ is $C$-tight in $\mathbb D([0,1])$.

Finally, deterministic scalar multiplication is continuous on $\mathbb D([0,1])$, so the rescaled process
\[
Z_t^T=\sqrt{\frac{T^\alpha \nu^*\delta^{-1}}{1-a_T}}\,(Y_t^T-\Lambda_t^T)
\]
is also $C$-tight.

We now prove tightness of $(Z^T)$.
Recall that
\[
Z_t^T=\sqrt{\frac{T^\alpha \nu^*\delta^{-1}}{1-a_T}}\,(Y_t^T-\Lambda_t^T)
      =\sqrt{\frac{1-a_T}{T^\alpha \nu^*\delta^{-1}}}\; M_{\lfloor tT\rfloor}^T,
\]
where $M_n^T=\sum_{s=1}^n(X_s^T-\lambda_s^T)$ is the compensated martingale.
Set
\[
D_T:=\frac{1-a_T}{T^\alpha \nu^*\delta^{-1}}.
\]
Then the jump at time $s/T$ is
\[
\Delta Z_{s/T}^T=\sqrt{D_T}\,(X_s^T-\lambda_s^T).
\]
Using that $X_s^T\mid\mathcal F_{s-1}\sim \mathrm{Pois}(\lambda_s^T)$, the conditional fourth centered moment satisfies
\[
\mathbb E\!\left[(X_s^T-\lambda_s^T)^4\mid\mathcal F_{s-1}\right]
= \lambda_s^T+3(\lambda_s^T)^2,
\]
hence, by the moment bounds on $\lambda_s^T$ proved in Lemma~\ref{lemma:tightnessofyandlambda},
\[
\sup_{1\le s\le \lfloor T\rfloor}\mathbb E\!\left[(X_s^T-\lambda_s^T)^4\right]
\le C\,T^{4\alpha-2}
\qquad(T\ge1).
\]
Therefore, by a union bound and Markov's inequality,
\begin{align*}
\mathbb P\!\left(\sup_{t\in[0,1]}|\Delta Z_t^T|>\varepsilon\right)
\le &\sum_{s=1}^{\lfloor T\rfloor}\frac{\mathbb E[(\Delta Z_{s/T}^T)^4]}{\varepsilon^4}\\
\le &\frac{T\,D_T^2}{\varepsilon^4}\,\sup_{s\le \lfloor T\rfloor}\mathbb E[(X_s^T-\lambda_s^T)^4]\\
\le &\frac{C}{\varepsilon^4}\,T\cdot T^{-4\alpha}\cdot T^{4\alpha-2}
= \frac{C}{\varepsilon^4}\,T^{-1}\to 0.
\end{align*}
In particular, $\sup_{t\in[0,1]}|\Delta Z_t^T|\xrightarrow{\mathbb P}0$.

Next, the predictable quadratic variation of $Z^T$ is
\[
\langle Z^T\rangle_t
=\sum_{s=1}^{\lfloor tT\rfloor}\mathbb E\!\left[(\Delta Z_{s/T}^T)^2\mid\mathcal F_{s-1}\right]
= D_T\sum_{s=1}^{\lfloor tT\rfloor}\lambda_s^T
=\Lambda_t^T.
\]
Since $(\Lambda^T)$ is $C$-tight by Lemma~\ref{lemma:tightnessofyandlambda}, the family
$(\langle Z^T\rangle)$ is $C$-tight in $\mathbb D([0,1])$ and any limit point is a.s.\ continuous.

Moreover, $\sup_{t\in[0,1]}|\Delta Z_t^T|\xrightarrow{\mathbb P}0$, with $\Delta Z_t^T:=Z_t^T-Z_{t-}^T$,
i.e.\ the ARJ$(1)$ condition of \cite{rebolledo1980central}. Therefore, by \cite[Proposition~12]{rebolledo1980central},
$(Z^T)$ is $C$-tight in $\mathbb D([0,1])$.

Finally, since both $(Z^T)$ and $(Y^T)$ are tight in $\mathbb D([0,1])$, the pair $(Z^T,Y^T)$ is tight in $\mathbb D([0,1])^2$.
}

\textbf{2. Convergence of the Quadratic Variation.}
Let $(Z,Y)$ be any limit point of $(Z^T, Y^T)$. To prove that $Y=[Z,Z]$, we follow the standard procedure of showing that both $[Z^T, Z^T]$ and $Y^T$ have the same limit. {More precisely, we prove the following two convergences in probability, uniformly on $[0,1]$:}
\[
  {\sup_{t\in[0,1]}\left|[Z^T, Z^T]_t - \langle Z^T, Z^T \rangle_t\right| \xrightarrow[T\to\infty]{\mathbb P} 0,}
\]
and
\[
  {\sup_{t\in[0,1]}\left|\langle Z^T, Z^T \rangle_t - Y^T_t\right| \xrightarrow[T\to\infty]{\mathbb P} 0.}
\]
{In particular, for each fixed $t\in[0,1]$, these convergences hold in probability at time $t$.}

Let $C_T^2 = \frac{T^\alpha \nu^* \delta^{-1}}{1-a_T}$ and $D_T = \frac{1-a_T}{T^\alpha \nu^* \delta^{-1}}$. Note that 
\[
  Z_t^T = C_T (Y_t^T - \Lambda_t^T) = C_T D_T M_{\lfloor tT \rfloor}^T.
\]
The increment is $\Delta Z_{s/T}^T = C_T D_T (X_s^T - \lambda_s^T)$.

{
(a) \textit{Asymptotic equivalence of $[Z^T, Z^T]$ and $\langle Z^T, Z^T \rangle$.}
}
{
Note that both $[Z^T,Z^T]$ and $\langle Z^T,Z^T\rangle$ depend on $T$; what we prove is that their difference is asymptotically negligible:
\[
\sup_{t\in[0,1]}\left|[Z^T,Z^T]_t-\langle Z^T,Z^T\rangle_t\right|\xrightarrow[T\to\infty]{\mathbb P}0.
\]
}
{
The predictable quadratic variation of $Z^T$ is given by the sum of conditional second moments of its increments.
\begin{align*}
  \langle Z^T, Z^T \rangle_t 
  = &\sum_{s=1}^{\lfloor tT \rfloor} \mathbb{E}\!\left[\left.(\Delta Z_{s/T}^T)^2 \,\right|\, \mathcal{F}_{s-1}\right] \\
  = &(C_T D_T)^2 \sum_{s=1}^{\lfloor tT \rfloor} \mathbb{E}\!\left[\left.(X_s^T - \lambda_s^T)^2 \,\right|\, \mathcal{F}_{s-1}\right] \\
  = &(C_T D_T)^2 \sum_{s=1}^{\lfloor tT \rfloor} \mathrm{Var}\!\left(\left.X_s^T \,\right|\, \mathcal{F}_{s-1}\right) \\
  = &(C_T D_T)^2 \sum_{s=1}^{\lfloor tT \rfloor} \lambda_s^T,
\end{align*}
where we used that, by construction, $X_s^T\mid \mathcal{F}_{s-1}\sim \mathrm{Pois}(\lambda_s^T)$, hence $\mathrm{Var}(X_s^T\mid \mathcal{F}_{s-1})=\lambda_s^T$.}
The sharp bracket (quadratic variation) is 
\[
  [Z^T, Z^T]_t = \sum_{s=1}^{\lfloor tT \rfloor} (\Delta Z_{s/T}^T)^2 = (C_T D_T)^2 \sum_{s=1}^{\lfloor tT \rfloor} (X_s^T - \lambda_s^T)^2.
\]
Consider the process 
\[
  U_n^T := [Z^T, Z^T]_{n/T} - \langle Z^T, Z^T \rangle_{n/T}.
\]
This is a martingale. We show it converges to zero in $L^2$:
\begin{align*}
    \mathbb{E}[(U_{\lfloor tT \rfloor}^T)^2] &= (C_T D_T)^4 \mathbb{E}\left[ \left( \sum_{s=1}^{\lfloor tT \rfloor} \left( (X_s^T - \lambda_s^T)^2 - \lambda_s^T \right) \right)^2 \right] \\
    &= (C_T D_T)^4 \sum_{s=1}^{\lfloor tT \rfloor} \mathbb{E}\left[ \left( (X_s^T - \lambda_s^T)^2 - \lambda_s^T \right)^2 \right] \\
    &= (C_T D_T)^4 \sum_{s=1}^{\lfloor tT \rfloor} \mathbb{E}[\text{Var}((X_s^T-\lambda_s^T)^2|\mathcal{F}_{s-1})] \\
    &= (C_T D_T)^4 \sum_{s=1}^{\lfloor tT \rfloor} \mathbb{E}[\lambda_s^T + 2(\lambda_s^T)^2].
\end{align*}
To show this convergence, we analyze the expectation of the sum. From the moment calculations in the proof of Lemma \ref{lemma:tightnessofyandlambda}, we know there exist uniform constants $\bar C_1,\bar C_2$ such that for all $s$ and $T$:
\[ \mathbb{E}[\lambda_s^T] \le \bar{C}_1 T^{2\alpha-1} \quad \text{and} \quad \mathbb{E}[(\lambda_s^T)^2] \le ((\bar C_1)^2+\bar C_1\bar C_2) T^{4\alpha-2}. \]
Therefore, the sum can be bounded as follows:
\begin{align*}
    \sum_{s=1}^{\lfloor tT \rfloor} \mathbb{E}[\lambda_s^T + 2(\lambda_s^T)^2] 
    &\le \sum_{s=1}^{\lfloor tT \rfloor} (\bar{C}_1 T^{2\alpha-1} + 2((\bar C_1)^2+\bar C_1\bar C_2) T^{4\alpha-2}) \\
    &\le \lfloor tT \rfloor \cdot (\bar{C}_1 T^{2\alpha-1} + 2((\bar C_1)^2+\bar C_1\bar C_2) T^{4\alpha-2}) \\
    &\le (\bar{C}_1 + 2((\bar C_1)^2+\bar C_1\bar C_2)) T^{4\alpha-1} \quad (\text{since } 2\alpha-1 < 4\alpha-2\ \text{and}\ t\in[0,1]).
\end{align*}

For the pre-factor, by Assumption \ref{assumption:2}, there exists a constant $\bar C_6$ such that for all sufficiently large $T$:
\[ (C_T D_T)^4 = D_T^2 = \left( \frac{1-a_T}{T^\alpha \nu^*\delta^{-1}} \right)^2 \le \bar C_6 T^{-4\alpha}. \]

Combining these non-asymptotic bounds, we have:
\begin{align*}
    \mathbb{E}[(U_{\lfloor tT \rfloor}^T)^2] 
    &\le (\bar C_6 T^{-4\alpha}) \cdot \left( (\bar{C}_1 + 2((\bar C_1)^2+\bar C_1\bar C_2)) T^{4\alpha-1} \right) \\
    &= \bar C_6 (\bar{C}_1 + 2((\bar C_1)^2+\bar C_1\bar C_2)) T^{-1}.
\end{align*}
As $T \to \infty$, this upper bound for the $L^2$-norm converges to 0. Convergence in $L^2$ implies convergence in probability, thus $[Z^T, Z^T]_t - \langle Z^T, Z^T \rangle_t \to 0$ in probability, uniformly in $t$ by Doob's inequality.

{
(b) \textit{Asymptotic equivalence of $\langle Z^T, Z^T \rangle$ and $Y^T$ (same limit points).}
}
{
Our goal is to show that the difference between $\langle Z^T, Z^T \rangle$ and $Y^T$ vanishes in probability, uniformly on $[0,1]$:
\[
\sup_{t\in[0,1]}\left|\langle Z^T, Z^T \rangle_t-Y_t^T\right|\xrightarrow[T\to\infty]{\mathbb P}0.
\]
Equivalently, $\langle Z^T, Z^T \rangle$ and $Y^T$ have the same limit points in $\mathbb D([0,1])$.
}

First, we establish a crucial identity. From the definitions of $C_T$ and $D_T$, we have
\[
    (C_T D_T)^2
    =\left(\frac{T^\alpha \nu^*\delta^{-1}}{1-a_T}\right)\left(\frac{1-a_T}{T^\alpha \nu^*\delta^{-1}}\right)^2
    =\frac{1-a_T}{T^\alpha \nu^*\delta^{-1}}
    =D_T.
\]
Hence the predictable quadratic variation is
\[
\langle Z^T, Z^T \rangle_t
=\sum_{s=1}^{\lfloor tT\rfloor}\mathbb E\!\left[(\Delta Z_{s/T}^T)^2\mid\mathcal F_{s-1}\right]
=(C_TD_T)^2\sum_{s=1}^{\lfloor tT\rfloor}\lambda_s^T
=D_T\sum_{s=1}^{\lfloor tT\rfloor}\lambda_s^T
=\Lambda_t^T.
\]
Therefore,
\[
\sup_{t\in[0,1]}\left|\langle Z^T, Z^T \rangle_t-Y_t^T\right|
=\sup_{t\in[0,1]}\left|\Lambda_t^T-Y_t^T\right|.
\]
By Lemma~\ref{lemma:YLambdaconvergeinprobability}, we already have
\[
\sup_{t\in[0,1]}\left|\Lambda_t^T-Y_t^T\right|\xrightarrow[T\to\infty]{\mathbb P}0,
\]
and thus
\[
\sup_{t\in[0,1]}\left|\langle Z^T, Z^T \rangle_t-Y_t^T\right|\xrightarrow[T\to\infty]{\mathbb P}0,
\]
which proves the claimed asymptotic equivalence (and in particular that $\langle Z^T\rangle$ and $Y^T$
share the same limit points).

\textbf{3. Characterization of the Limit.}
Since $(Z^T, Y^T)$ converges in distribution to $(Z,Y)$, and $[Z^T, Z^T]_t \to Y_t$ in probability, it follows from the properties of weak convergence that the limit process satisfies $Y=[Z,Z]$.

Furthermore, we have established that $(Z^T)$ is $C$-tight. A sequence of martingales that is $C$-tight and whose quadratic variation converges to a continuous process must converge to a continuous martingale (see, e.g., Proposition VI-3.26 in \cite{jacod2013limit}). Since $[Z^T,Z^T]_t \to Y_t$ and $Y$ (as the limit of a $C$-tight sequence $Y^T$) is continuous, the limit process $Z$ must be a continuous martingale.


\subsection{Proof of Theorem \ref{theorem:1}}
\label{subsec:proofoftheorem1}

In fact, we only need to verify that the limit process satisfies equation \eqref{eq:theorem1equation}. The proof of the H\"older property for $Y$ in Theorem \ref{theorem:1} is exactly the same as in \cite{jaisson2016rough}. {Let $(Z^{T_n}, Y^{T_n})$ be a subsequence converging in distribution to $(Z,Y)$, and for simplicity write $(Z^T,Y^T)$ instead of $(Z^{T_n},Y^{T_n})$. By Skorohod's representation theorem, there exists a probability space on which one can define copies (with the same laws) of $(Z^T,Y^T)$ and $(Z,Y)$ such that $(Z^T,Y^T)\to(Z,Y)$ almost surely in the Skorohod topology. We work on this auxiliary probability space from now on. In particular, to apply the Dambis--Dubins--Schwarz theorem below, we may pass to the usual extension of this probability space so that there exists a Brownian motion $B$ such that $Z_t = B_{Y_t}$.}
Since the processes $Z$ and $Y$ are continuous, we obtain 
\[
  \sup_{t\in[0,1]}|Y_t^T-Y_t|\rightarrow 0,\ \sup_{t\in[0,1]}|Z_t^T-Z_t|\rightarrow 0.
\]

  Let us now rewrite $\sum_{s=1}^n\lambda_s^T$. We can use Lemma \ref{lemma:discreterenewalequation} to obtain the expression of $\sum_{s=1}^n\lambda_s^T$, for all $n\in \mathbb N$, using Abel's lemma, we obtain
  \begin{align*}
    \sum_{s=1}^n\lambda_s^T
    =&n\mu^T+\sum_{s=1}^n\eta_{n-s}^T\sum_{u=1}^{s}\lambda_u^T+\sum_{s=1}^n\sum_{u=1}^{s-1}\eta_{s-u}^T(M_u^T-M_{u-1}^T)\\
    =&n\mu^T+\sum_{s=1}^n\eta_{n-s}^T\sum_{u=1}^s\lambda_u^T+\sum_{s=1}^n\eta_{n-s}^TM_s^T.
  \end{align*}
  
  Then we try to use Lemma \ref{lemma:discreterenewalequation}, we first calculate
  \begin{align*}
    \sum_{s=1}^{n-1}A_{n-s}^T\sum_{u=1}^s\eta_{s-u}^TM_u^T
    =&\sum_{u=1}^{n-1}M_u^T\sum_{s=1}^{n-u-1}A_{n-s-u}^T\eta_s^T
    =\sum_{u=1}^{n-1}M_u^T(A_{n-u}^T-\eta_{n-u}^T),
  \end{align*}
  where we use the fact that $A^T*\eta^T=A^T-\eta^T$. By Lemma \ref{lemma:discreterenewalequation}, we obtain
  \[
    \sum_{s=1}^n\lambda_s^T=n\mu^T+\mu^T\sum_{s=1}^{n-1}sA_{n-s}^T+\sum_{s=1}^{n-1}A_{n-s}^TM_s^T.
  \]
  Substituting $n$ with $\lfloor tT \rfloor$ and multiplying by $(1-a_T)/(T^\alpha \nu^*\delta^{-1})$, denote
  \[
    u_T=\frac{\mu^T}{\nu^*\delta^{-1}T^{\alpha-1}},
  \]
  we can then decompose $\Lambda_t^T$ into three parts
  \[
    \Lambda_t^T=\frac{1-a_T}{T^\alpha \nu^*\delta^{-1}}\sum_{s=1}^{\lfloor tT \rfloor}\lambda_s^T=T_1+T_2+T_3,
  \]
  where
  \begin{align*}
    &T_1=(1-a_T)tu_T,\\
    &T_2=\frac{1-a_T}{T^\alpha \nu^*\delta^{-1}}\mu^T\sum_{s=1}^{\lfloor tT \rfloor-1}sA_{\lfloor tT \rfloor-s}^T
    =a_Tu_T\sum_{s=1}^{\lfloor tT \rfloor-1}\frac1T\left( \frac{s}{T} \right)\rho^T\left(t-\frac{s}{T}\right),\\
    &T_3=\frac{1-a_T}{T^\alpha \nu^*\delta^{-1}}\sum_{s=1}^{\lfloor tT \rfloor-1}A_{\lfloor tT \rfloor-s}^TM_s^T,
  \end{align*}
  and $\rho^T(\cdot)$ is defined in \eqref{eq:rho_T_def}. {As $T \to \infty$, since $u_T \to 1$ and $(1-a_T) \to 0$, we have $T_1 \to 0$. 
For $T_2$, recognized as a Riemann sum converging to a convolution integral, and by integrating by parts:
\begin{align*}
  T_2 &\sim a_Tu_T\int_0^t\rho^T(t-s)sds = a_Tu_T\int_0^t F^T(t-s)ds \\
      &\rightarrow \int_0^t F^{\alpha,\lambda}(t-s)ds = \int_0^t f^{\alpha,\lambda}(t-s)s \, ds.
\end{align*}
Finally, we turn to $T_3$. Remark that $M_s^T = \sqrt{\frac{T^\alpha \nu^*\delta^{-1}}{1-a_T}}Z_{s/T}^T$. Thus,
\begin{align*}
  T_3 &= \frac{a_T}{\sqrt{T^\alpha (1-a_T)\nu^*\delta^{-1}}} \sum_{s=1}^{\lfloor tT \rfloor-1}\frac1{T}\rho^T\left(t-\frac{s}{T}\right)Z_{s/T}^T \\
      &\sim \frac{a_T}{\sqrt{T^\alpha (1-a_T)\nu^*\delta^{-1}}} \int_0^t\rho^T(t-s)Z_s^T ds.
\end{align*}
Using Assumption \ref{assumption:2}, the prefactor converges to:
\[
    \lim_{T\to\infty} \frac{1}{\sqrt{T^\alpha (1-a_T)\nu^*\delta^{-1}}} = \frac{1}{\sqrt{\lambda\delta \cdot \nu^*\delta^{-1}}} = \frac{1}{\sqrt{\lambda\nu^*}}.
\]
Following the convergence arguments in Theorem 3.1 of \cite{jaisson2016rough}, $T_3$ converges to
\[
  \frac1{\sqrt{\nu^*\lambda}}\int_0^t \lambda(t-s)^{\alpha-1}E_{\alpha,\alpha}(-\lambda(t-s)^\alpha)Z_s ds.
\]
Since $Z$ is a continuous martingale with quadratic variation $Y$, the fact that $Z_t=B_{Y_t}$ is a consequence of the Dambis-Dubin-Schwarz theorem (e.g., Theorem V-1.65 in \cite{revuz2013continuous}).}


\subsection{Proof of Corollary \ref{cor:initial_value}}
\label{subsec:proofofcorollaryinitialvalue}
The proof extends the arguments used to establish Theorem \ref{theorem:1}. The core of the argument lies in analyzing the convergence of the scaled process $Y_t^T$ under the new time-dependent baseline intensity $\hat{\mu}_T(n)$.

The proof proceeds by decomposing the process $Y_t^T$ into its deterministic mean, $\mathbb{E}[Y_t^T]$, and its centered stochastic part, $Y_t^T - \mathbb{E}[Y_t^T]$, and then analyzing the limit of each component separately.

Let us denote the deterministic mean as $\bar{\Lambda}_t^T := \mathbb{E}[Y_t^T]$. Using the expression for $\mathbb{E}[N_n^T]$ from Lemma \ref{lemma:renewal_expansion_with_initial_value} and applying the scaling factor, we can decompose $\bar{\Lambda}_t^T$ into several terms. After substituting the definition of $\hat{\mu}_T(n)$, the expression for $\mathbb{E}[N_n^T]$ becomes:
\[
\mathbb{E}[N_n^T] = \mu^T n + (1-\xi)\mu^T \sum_{s=1}^{n-1}sA_{n-s}^T + \xi\mu^T \left( \frac{n}{1-a_T} - \sum_{s=1}^{n-1}sA_{n-s}^T \right).
\]
Scaling this expression gives the components of the deterministic mean $\bar{\Lambda}_t^T$:
\begin{itemize}
    \item A term that vanishes in the limit:
        \[
            T_1 = \frac{(1-a_T)\mu^T \lfloor Tt \rfloor}{T^\alpha \nu^*\delta^{-1}}.
        \]
    \item The main deterministic term that converges to a non-trivial limit:
        \[
            T_2 = \frac{(1-a_T)\mu^T}{T^\alpha \nu^*\delta^{-1}} \left[ (1-\xi) \sum_{s=1}^{\lfloor Tt \rfloor-1}sA_{\lfloor Tt \rfloor-s}^T + \xi \left( \frac{\lfloor Tt \rfloor}{1-a_T} - \sum_{s=1}^{\lfloor Tt \rfloor-1}sA_{\lfloor Tt \rfloor-s}^T \right) \right].
        \]
\end{itemize}
The centered stochastic part, $Y_t^T - \bar{\Lambda}_t^T$, is obtained by scaling the expression from Equation \eqref{eq:centered_N_initial_value} in Lemma \ref{lemma:renewal_expansion_with_initial_value}. We denote this scaled stochastic part as $T_3^{\text{stoch}}$:
\begin{itemize}
    \item The scaled centered stochastic part:
        \[
            T_3^{\text{stoch}} := Y_t^T - \bar{\Lambda}_t^T = \frac{1-a_T}{T^\alpha \nu^*\delta^{-1}} \left( M_{\lfloor Tt \rfloor}^T + \sum_{s=1}^{\lfloor Tt \rfloor-1} A_{\lfloor Tt \rfloor-s}^T M_s^T \right).
        \]
\end{itemize}
Next, we analyze the limits of these terms as $T \to \infty$:
\begin{itemize}
    \item By Assumption \ref{assumption:2}, it is straightforward to show that $\lim_{T\to\infty} T_1 = 0$.
    \item For $T_2$, we use the same Riemann sum approximation argument as in the proof of Theorem \ref{theorem:1}, where the sum $\sum_{s=1}^{n-1} s A_{n-s}^T$ corresponds to the integral of the CDF. In the limit, the sum of the deterministic parts converges:
    \[
    \lim_{T\to\infty} \bar{\Lambda}_t^T = \lim_{T\to\infty} (T_1 + T_2) = (1-\xi)\int_0^t F^{\alpha,\lambda}(u)du + \xi \left( t - \int_0^t F^{\alpha,\lambda}(u)du \right).
    \]
    \item The convergence of the scaled stochastic part, $T_3^{\text{stoch}}$, was established in the proof of Theorem \ref{theorem:1}, as the structure of the martingale part is identical. Its limit is given by:
    \[
    \lim_{T\to\infty} T_3^{\text{stoch}} = \frac{1}{\sqrt{\nu^* \lambda}} \int_0^t f^{\alpha,\lambda}(t-s) B_{Y_s} ds.
    \]
\end{itemize}
Combining these limits, the limit process $Y_t = \lim_{T\to\infty} Y_t^T$ is given by the sum of the limits of its deterministic and stochastic parts:
\[
Y_t = (1-\xi)\int_0^t F^{\alpha,\lambda}(u)du + \xi \left( t - \int_0^t F^{\alpha,\lambda}(u)du \right) + \frac{1}{\sqrt{\nu^* \lambda}} \int_0^t f^{\alpha,\lambda}(t-s) B_{Y_s} ds.
\]
The process $\dot{Y}_t$ is the derivative of $Y_t$. Differentiating the deterministic part with respect to $t$ gives:
\[
\frac{d}{dt} \left( \lim_{T\to\infty} \bar{\Lambda}_t^T \right) = (1-\xi)F^{\alpha,\lambda}(t) + \xi\left(1 - F^{\alpha,\lambda}(t)\right) = F^{\alpha,\lambda}(t) + \xi\left(1 - F^{\alpha,\lambda}(t)\right).
\]
The derivative of the stochastic part, as established in Corollary \ref{corollary:derivativesde}, yields the stochastic integral term. Summing these derivatives gives the full stochastic Volterra equation for $\dot{Y}_t$ as stated in \eqref{eq:sve_with_initial}.

Finally, to confirm the initial condition, we evaluate at $t=0$. Since $F^{\alpha,\lambda}(0)=0$ and the stochastic integral is zero at $t=0$, we have:
\[
\dot{Y}_0 = F^{\alpha,\lambda}(0) + \xi\left(1 - F^{\alpha,\lambda}(0)\right) + 0 = \xi.
\]
This concludes the proof.

\end{document}